\documentclass[10pt]{amsart}
\usepackage{enumerate,amsmath, amsfonts, amssymb, latexsym, amscd, amsthm}

\newlength{\tabwidth}
\newlength{\tabheight}
\newlength{\tabrule}
\newlength{\tabwidthx}
\newlength{\tabheightx}

\def\gentabbox#1#2#3#4{\vbox to \tabheight{\setlength{\tabrule}{#3}%
  \setlength{\tabwidthx}{#1\tabwidth}\addtolength{\tabwidthx}{\tabrule}%

\setlength{\tabheightx}{#2\tabheight}\addtolength{\tabheightx}{-\tabheight}%
  \hbox to #1\tabwidth{%
    \hspace{-0.5\tabrule}\rule{\tabrule}{#2\tabheight}\hspace{-\tabrule}%
    \vbox to #2\tabheight{\hsize=\tabwidthx%
      \vspace{-0.5\tabrule}\hrule width\tabwidthx height\tabrule%
      \vspace{-0.5\tabrule}\vfil%
      \hbox to \tabwidthx{\hss#4\hss}%
        \vfil\vspace{-0.5\tabrule}%
      \hrule width\tabwidthx height\tabrule\vspace{-0.5\tabrule}}%
    \hspace{-\tabrule}\rule{\tabrule}{#2\tabheight}\hspace{-0.5\tabrule}}%
  \vspace{-\tabheightx}}}
\def\genblankbox#1#2{\vbox to \tabheight{\vfil\hbox to
#1\tabwidth{\hfil}}}
\def\tabbox#1#2#3{\gentabbox{#1}{#2}{0.4pt}{\strut #3}}

\catcode`\:=13 \catcode`\.=13 \catcode`\;=13 \catcode`\>=13
\catcode`\^=13
\def:#1\\{\hbox{$#1$}}
\def.#1{\tabbox{1}{1}{$#1$}}
\def>#1{\tabbox{2}{1}{$#1$}}
\def^#1{\tabbox{1}{2}{$#1$}}
\def;{\genblankbox{1}{1}\relax}
\catcode`\:=12 \catcode`\.=12 \catcode`\;=12 \catcode`\>=12
\catcode`\^=7

\newenvironment{tableau}{\bgroup\catcode`\:=13 \catcode`\.=13
  \catcode`\;=13 \catcode`\>=13 \catcode`\^=13
  \setlength{\tabheight}{3ex}\setlength{\tabwidth}{3ex}%
  \def\b##1##2##3{\gentabbox{##1}{##2}{1.2pt}{\vbox{##3}}}%
  \def\n##1##2##3{\gentabbox{##1}{##2}{0.4pt}{\vbox{##3}}}%
  \vbox\bgroup\offinterlineskip}{\egroup\egroup}

\input xy
\xyoption{all}

\newtheorem{fact}{Fact}[section]
\newtheorem{theorem}[fact]{Theorem}
\newtheorem{corollary}[fact]{Corollary}
\newtheorem{lemma}[fact]{Lemma}
\newtheorem*{theorem*}{Theorem}
\newtheorem*{fact*}{Fact}
\newtheorem*{corollary*}{Corollary}
\newtheorem*{lemma*}{Lemma}
\newtheorem{proposition}[fact]{Proposition}

\theoremstyle{definition}
\newtheorem{definition}[fact]{Definition}
\newtheorem*{definition*}{Definition}

\theoremstyle{remark}

\newtheorem{example}[fact]{Example}

\newtheorem*{remark*}{Remark}
\newtheorem*{exercise*}{Exercise}
\newtheorem*{question*}{Question}
\newtheorem*{example*}{Example}
\newtheorem*{problem*}{Problem}

\newcommand\lieg{\mathfrak {g}}

\newcommand\calV{\mathcal{V}}
\newcommand\calO{\mathcal{O}}
\newcommand\calU{\mathcal{U}}
\newcommand{\V}{\vspace{0.4cm}}

\newcommand{\calC}{\mathcal{C}}

\newcommand{\Irr}{\text{Irr}}

\title{Orbital Varieties and Unipotent Representations}
\author{Thomas Pietraho}
\date{}

\begin{document}

 \maketitle

\begin{abstract}
Using the notion of a Lagrangian covering, W.~Graham and D.~Vogan proposed a method of constructing representations from the coadjoint orbits for a complex semisimple Lie group $G$.  When the coadjoint orbit $\calO$ is nilpotent,  a representation of $G$ is attached to each orbital variety of $\calO$ in this way.  In the setting of classical groups, we show that whenever it is possible to carry out the Graham-Vogan
construction for an orbital variety of a spherical $\mathcal{O}$, its infinitesimal character lies in a set of characters attached to $\mathcal{O}$ by W.~M.~McGovern. Furthermore, we show
that it is possible to carry out the Graham-Vogan construction for a
sufficient number of orbital varieties to account for all the
infinitesimal characters in this set.
\end{abstract}

\section{Introduction}

Consider a complex semisimple Lie group $G$.  The philosophy of the orbit method seeks to parameterize the unitary dual of $G$ via  orbits of the coadjoint action of $G$ on the dual of its Lie algebra.  The representations attached to nilpotent coadjoint orbits are the so-called unipotent representations of $G$.  In \cite{gv}, W.~Graham and D.~Vogan have proposed a general method of attaching a set of representations to each coadjoint orbit $\calO$ of $G$.  Very little is known about which representations of $G$  actually arise in this way, but conjecturally, when $\calO$ is nilpotent, they should coincide with the set of unipotent representations corresponding to  $\calO$.
The goal of this paper is to shed some light on the situation, examining the Graham-Vogan construction for the family of spherical nilpotent orbits of $G$.

The idea behind the construction of \cite{gv} is to generalize the method of polarization.  Instead of using Lagrangian foliations, it relies on the more general notion of a Lagrangian covering of a coadjoint orbit.
When $\calO$ is nilpotent, the main ingredients  of a Lagrangian covering are certain Lagrangian submanifolds called {\it orbital varieties}.  For
each choice of orbital variety $\calV$ and a choice of an {\it
admissible orbit datum} $\pi$, the  Graham-Vogan construction
defines a subspace $V(\calV, \pi)$ of sections of a bundle
over $G/Q_{\calV}$, where $Q_{\calV}
\subset G$ is the maximal subgroup of $G$ stabilizing $\calV$.

To examine which representations of $G$ arise among the spaces $V(\calV, \pi)$, we rely on a combinatorial description of orbital varieties in classical groups obtained by W.~M.~McGovern, \cite{mcgovern2}, as well as the author, \cite{pietraho:springer}.  They are parametrized by standard Young tableaux in type $A$ and standard domino tableaux in the other classical types.  There are two advantages to this description.  First, it is  easy to determine $Q_{\calV}$ from the tableau parameterizing $\calV$.  Second, the parametrization itself suggests a means of addressing our work inductively, setting up a framework for our calculations.

As our main goal is to determine how well the $V(\calV, \pi)$ fit the role of unipotent representations, we would like to describe what reasonable conditions for this might be.  There is at least one commonly accepted necessary criterion for a representation $V$ to be attached to an orbit $\calO$.  According to work of Borho and Brylinski, the variety $\mathcal{V}(\textup{Ann}(V)) \subset \mathfrak{g}^*$  associated to the annihilator of $V$ in the universal enveloping algebra of $\mathfrak{g}$ is the closure of a single nilpotent orbit if $V$ is irreducible.
Thus, for a unipotent representation $V$ arising from the nilpotent orbit $\calO$, we should expect
$$\mathcal{V}(\textup{Ann}(V)) = \overline{\calO}.$$

A classification of unitary representations of complex reductive Lie groups can be obtained from a construction which begins with a set of {\it special} unipotent representations first suggested by J.~Arthur, see \cite{barbasch}.  However, only {\it special} nilpotent orbits arise as associated varieties of  special unipotent representations.  To remedy this shortfall, McGovern has  suggested extending the set of special unipotent representations to a set of $q$-unipotent representations whose associated varieties include all nilpotent orbits of $G$ \cite{mcgovern}.  Included in his work is a description of the infinitesimal characters of $q$-unipotent representations for classical groups, suggesting a natural benchmark for examining the Graham-Vogan spaces.
After incorporating certain geometric considerations into McGovern's list, we will define a set $IC^1({\calO})$ of infinitesimal characters attached to each nilpotent orbit $\mathcal{O}$.  The set of infinitesimal characters of
the representations attached to $\calO$ should contain $IC^1({\calO})$.  Our main result is that for spherical $\calO$, this is exactly what happens.  We paraphrase this as follows:

\vspace{.1in}
\noindent
{\bf Theorems \ref{theorem:all} and \ref{theorem:ic}.} {\it Let
$\calO$ be a spherical nilpotent orbit of a complex classical
semisimple Lie group $G$ of rank $n$ and write $GV_{\calO}$ for the set of representations $V(\calV, \pi)$ arising for some $\mathcal{V} \subset \calO$. Let $\chi_{\calV}$ be the
infinitesimal character associated to $V(\calV, \pi)$. Then,
      \begin{enumerate}[(i)]
        \item If $\calO$ is rigid, then  $IC^1({\calO}) = \{\chi_{\calV}
             \; | \; V(\calV, \pi ) \in GV_{\calO}\},$

        \item If $\calO$ is a model orbit and $n>2$, then   $IC^1({\calO})
         \subset \{\chi_{\calV}|\; V(\calV, \pi) \in GV_{\calO}\}.$
      \end{enumerate}
}
\vspace{.1in}

The above theorems imply that, at least for spherical nilpotent orbits,
the Graham-Vogan spaces are indeed candidates for unipotent
representations.  For larger non-spherical orbits, calculations indicate that the family of infinitesimal characters of representations in $GV_{\calO}$ is too numerous to
form the set of representations attached to $\mathcal{O}$.  However, additional conditions on the
closure $\overline{\calO}$ not considered in the construction \cite{gv} should make
it possible to restrict the resulting set of possible infinitesimal characters.

The paper is structured as follows.  Section \ref{section:construction} presents a summary of the construction of the spaces $V(\calV, \pi)$ and details the combinatorics of nilpotent coadjoint orbits and orbital varieties in classical groups.  Section \ref{section:restrictiontospherical} begins with an example which the rest of the paper is designed to mimic.  A crucial assumption is that the stabilizing parabolic $Q_\mathcal{V}$ has a dense orbit in $\calV$.  We restrict our attention to nilpotent orbits all of whose orbital varieties enjoy this property and detail this restriction in terms of the combinatorics of Section  \ref{section:construction}.  The section concludes with a description of the inductive process we will use in the rest of the paper.  Finally, Section \ref{section:infinitesimal} addresses infinitesimal characters.  We begin by detailing conditions under which it is possible to carry out the Graham-Vogan construction.  After defining the desired set $IC^1(\mathcal{O})$ of infinitesimal characters that ought to be attached to the orbit $\calO$, we compute which ones arise from spaces of the form $V(\calV,\pi)$.


\section{Preliminaries}

We begin this section with a brief outline of the Graham-Vogan construction.  In the setting of classical groups, both nilpotent orbits and orbital varieties, on which this construction relies, admit combinatorial descriptions. We summarize these and list a few useful results.

\subsection{The Graham-Vogan Construction}\label{section:construction}

Let $G$ be a semisimple Lie group and $\mathfrak{g}$ its Lie algebra.  The coadjoint orbit through a point $f \in \mathfrak{g}^*$ is the set $$\mathcal{O}_f = G \cdot f \cong G/G_f$$
where we write $G_f$ for the isotropy subgroup of the coadjoint action of $G$ based at $f$.  The nondegeneracy of the Killing form permits us to identify the set of coadjoint orbits with the set of adjoint orbits.  The Graham-Vogan construction of representations associated to a
coadjoint orbit $\calO$ is an extension of the method of polarizing
 a coadjoint orbit.  We briefly recount this work, following \cite{gv}.  It begins with the notion of a Lagrangian covering.

\begin{definition}
A {\it Lagrangian covering} of a symplectic manifold $\calO$ is a
pair $(Z,M)$  of manifolds and smooth maps $(\tau, \rho)$
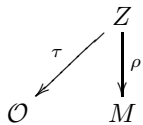
$$
\xymatrix{ & Z \ar[dl]_{\tau} \ar[d]^{\rho}  \\
           \calO  & M \\}
$$
such that the diagram is a double fibration and each fiber of $\rho$ is a Lagrangian submanifold of $\calO$.
\end{definition}

\begin{theorem}[\cite{ginsburg}]  Let $G$ be a complex reductive Lie group
  and $\calO$ be a coadjoint orbit.  Then there exists an equivariant
  Lagrangian covering of $\calO$ where $M$ is a partial flag variety for
  $G$.
\end{theorem}

We are interested in this construction when $\calO$ is a nilpotent coadjoint
orbit. Fix a Borel subgroup $B$ of $G$ with unipotent radical $N$.
Write
$\mathfrak{g}=\mathfrak{n}^- \oplus \mathfrak{t} \oplus
\mathfrak{n}$ for the corresponding triangular decomposition. Let
us restrict our attention to  nilpotent coadjoint orbits $\calO,$
and consider the set $\calO \cap \mathfrak{n}$. This is a locally
closed subset of $\mathfrak{n}$ and can be expressed as a union of
its irreducible components.

\begin{definition}
Consider a nilpotent coadjoint orbit $\calO.$ Denote the set of
irreducible components of the variety $\calO \cap \mathfrak{n}$ by
$\text{Irr}(\calO \cap \mathfrak{n})$.  Each element of
$\text{Irr}(\calO
  \cap \mathfrak{n})$ is an {\it orbital variety } for $\calO$.
\end{definition}

\begin{proposition}
The set $\text{Irr}(\calO \cap \mathfrak{n})$ is finite.  Further,  every orbital variety $\mathcal{V} \in \text{Irr}(\calO \cap \mathfrak{n})$ has $\dim \mathcal{V} = \frac{1}{2} \dim \mathcal{O}$ and
is a Lagrangian subvariety of $\calO$.
\end{proposition}
We will construct  a distinct Lagrangian covering for each
orbital variety contained in $\calO$. Fix an orbital variety $\calV$
and let $\calV^0$ be its smooth part. Let $Q= Q_{\calV}=  \{q \in G \, | \, q
\cdot \calV = \calV \}.$  This is a parabolic subgroup of $G$ since
$\calV$ is $B$-stable. Furthermore, define the manifold $M$
by $M=\{g \cdot \calV \, | \, g \in G \} \cong G/Q.$  It is a
partial flag variety for $G$.

\begin{definition}
 For a subgroup $H \subset G$ and an $H$-space
$V$, let $G \times_H V$ to be the set of equivalence classes
in $G \times V$ with $(gh, v) \sim (g, h \cdot v)$ for $g \in G, h
\in H,$ and $v \in V$.
\end{definition}

The manifold $Z$ in the Lagrangian covering of $\calO$ associated to
the orbital variety $\mathcal{V}$ is now defined to be  $Z= G \times_Q
\calV^0.$
The map $\rho: Z \longrightarrow M$ arises from the projection of $G$
onto $G/Q$. The action of $G$ on $\calO$ gives natural map $G\times
\calV \rightarrow \calO$.  It descends to an algebraic map $\tau: Z
\longrightarrow \calO$.  We now have a Lagrangian covering:
$$
\xymatrix{ & G \times_Q \calV^0 \ar[dl]_{\tau} \ar[d]^{\rho}  \\
            G/G_f &  G/Q \\
           }
$$
Because the diagram is a double fibration, we can identify fibers of
$\rho$ with subsets of $\calO$.  In fact, each fiber is Lagrangian
in $\calO$.

The next step is to construct a representation from this Lagrangian covering.  Suppose that we have a
$G$-equivariant line bundle $\mathcal{L}_M \rightarrow M$.  We can
again pull this bundle back along the fibration $\rho$, this time to
obtain a bundle $\mathcal{L}_Z$.

$$\xymatrix{  & G \times_Q \calV^0 \ar[d]^{\rho} \ar[dl]_{\tau}  & \mathcal{L}_Z \ar[l]\\
           \calO &  M & \mathcal{L}_M \ar[l] \ar[u]_{\rho^*} \\
           }
$$
Geometric quantization suggests that the representations attached to
$\calO$ should  lie in the space of sections of $\mathcal{L}_M$, or
in other words, in the space of sections of $\mathcal{L}_Z$ that are
constant on the fibers of $\rho.$  This is very similar to the
situation arising in the polarization construction, as the fibers of
$\rho$ can again be identified with Lagrangian submanifolds of
$\calO.$   As described in \cite{gv}, however, the full set of sections of
$\mathcal{L}_M$ is too large to quantize $\calO$ and we pick out a subspace.

We relate only a general overview, and
direct the reader to \cite{gv} itself for the relevant details.  The
main idea is to prune the full space of sections of $\mathcal{L}_M$,
leaving ones which also come from an {\it admissible orbit  datum}
of $\calO$.

To do this, one must first attach a geometric structure to each
orbit datum.  This is achieved by mimicking the construction of a
Hermitian bundle that often arises in descriptions of geometric
quantization of {\it
  integral} orbit data.  The main difficulty then lies in finding a way of
embedding the information from this bundle into the space of
sections of $\mathcal{L}_M$.

\begin{definition}
An {\it admissible orbit datum} at $f \in \mathfrak{g}^*$ is a
genuine irreducible unitary representation $\pi$ of the metaplectic
cover $\widetilde{G}_f$ satisfying $$\pi(\exp Y) = \chi(f(Y))$$ for
a fixed non-trivial character $\chi$ of $\mathbb{R}$.
\end{definition}

Denote the metaplectic representation of
$\widetilde{G}_f$ by $\tau_f$ and form the tensor product
representation $\pi \otimes \tau_f$.  While $\tau_f$ and $\pi$ are
genuine representations of $\widetilde{G}_f$, $\pi \otimes \tau_f$
in fact descends to a representation of $G_f$ itself.  This allows
us to define a Hilbert bundle over the coadjoint orbit $\calO$ by
$$\mathcal{S}_\pi = G \times_{G_f} (\pi \otimes \tau_f).$$ This is the {\it bundle of twisted symplectic spinors}
on $\calO$. The metaplectic representation $\tau_f$ of
$\widetilde{G}_f$ decomposes into two irreducible and inequivalent
representations $\tau_f^{odd}$ and $\tau_f^{even}.$  Write
$\tau_f^{odd, \infty}$ and $\tau_f^{even, \infty}$ for the
corresponding sets of smooth vectors. This decomposition passes to
the bundle $\mathcal{S}_\pi$ and the geometric structure attached to
the admissible orbit datum $\pi$ is the subbundle of
$\mathcal{S}_\pi$ defined by
$$\mathcal{S}_\pi^{even, \infty} = G \times_{G_f} (\pi \otimes \tau_f^{even,\infty}).$$

\begin{definition}
Suppose that $X$ is a symplectic manifold.  The {\it bundle of
  infinitesimal Lagrangians} on $X$ is a fiber bundle $\mathcal{B}(X)$
over $X$.  The fiber over each point $x \in X$ is the set of
Lagrangian subspaces of the tangent space at $x$ of $X$, denoted by
$\mathcal{B}(T_x{X})$.
\end{definition}

\begin{definition}
Let $\calO$ be a coadjoint orbit, and consider a Lagrangian $\calV$
in the tangent space $\mathfrak{g}/\mathfrak{g}_f$.  Write
$\mathcal{L}(\calV)$ for the line defined in \cite[7.4(c)]{gv} from
the metaplectic representation $\tau_f$.  The admissible orbit datum
$\pi$ defines a $G$-equivariant vector bundle $\calV_\pi$ on
$\mathcal{B}(\calO)$ by letting the fiber at each $\calV$ be
$\mathcal{H}_{\pi} \otimes \mathcal{L}(\calV).$
\end{definition}

\begin{theorem}[\cite{gv},\cite{kostant}]
There exists a natural inclusion
$$i: C^{\infty} (\calO, \mathcal{S}^{even, \infty}_\pi) \hookrightarrow
C^{\infty}(\mathcal{B}(\calO), \mathcal{V}_\pi).$$
\end{theorem}

Next, we incorporate the bundle $\calV_\pi$ over
$\mathcal{B}(\calO)$ into the Lagrangian covering diagram.  Define a
map $\sigma: Z \rightarrow \mathcal{B}(\calO)$ as follows.  Fix $z
\in Z$. The definition of Lagrangian covering forces the fiber of
$\rho$ over $\rho(z) \in M$ to be  a Lagrangian submanifold of
$\calO$ that contains $\tau(z)$.  Hence its tangent space
$T_{\tau(z)}(\rho^{-1}(\rho(z))$ is a Lagrangian subspace of
$T_{\tau(z)}(\calO)$ and thus an element of $\mathcal{B}(\calO)$.
Let $$ \sigma(z) = T_{\tau(z)}(\rho^{-1}(\rho(z)).$$  In this way,
$\sigma$ becomes a bundle map over $\calO$.

We can pull back the bundle
$\calV_\pi$ along $\sigma$ to a bundle $\sigma^*(\calV_\pi)$ over
$Z$.  Smooth sections of $\calV_\pi$ pull back to smooth sections of
$\sigma^*(\calV_\pi)$ and we have an injective map $ \sigma^* \cdot
i : C^\infty(\calO, \mathcal{S}_\pi^{even,
  \infty}) \hookrightarrow C^\infty (Z, \sigma^*(\calV_\pi)).$
{\it Provided that} there is a $G$-equivariant vector bundle isomorphism
$j_{\pi}: \sigma^*(\calV_\pi) \rightarrow \rho^*(\mathcal{L}_M)$ we
can define a smooth representation of $G$ as:
$$V(\calV, \pi)= \rho^*(C^\infty(M,\mathcal{L}_M))\cap
j_{\pi}(\sigma \cdot i(C^\infty(\calO, \mathcal{S}_\pi^{even,
\infty})))$$

If $\mathcal{L}_M$ is given by a representation $\gamma$ of the
parabolic subgroup $Q$, then $V(\calV, \pi)$ lies in the space of
smooth vectors of the degenerate principal series representation
induced from $\gamma.$  The entire construction may be summarized by
the following diagram.
$$
\xymatrix{ &\calV_\pi \ar[d] \ar[r]^{\sigma^*} & \sigma^*(\calV_\pi) \ar[d] \ar[dr]^{j_{\pi}} & \\
S_\pi^{even, \infty} \ar @/_/ [dr] & \mathcal{B}(\calO) \ar[d] &  G \times_Q \calV^0 \ar[l]_{\sigma} \ar[d]^{\rho} \ar[dl]_{\tau} & \mathcal{L}_Z \ar[l]\\
  &           \calO   &        G/Q & \mathcal{L}_M \ar[u]_{\rho^*} \ar[l] \\
            }
$$

\subsection{Nilpotent Orbits in Classical Types}

The nilpotent coadjoint orbits of classical complex simple Lie groups are parametrized by partitions.  When $G$ is of type $A$, Jordan block sizes determine nilpotent adjoint orbits leading to a one-to-one correspondence between them and partitions. In classical groups not of
type $A,$  the parameterization has the same flavor, but not all partitions arise from Jordan block decompositions.  We describe the details presently.

To be specific, let $\epsilon = \pm 1$ and take $<,>_{\epsilon}$ to be a non-degenerate
bilinear form on $\mathbb{C}^m$ such that $$<x,y>_{\epsilon} =
\epsilon <y,x>_{\epsilon} \forall x,y \in \mathbb{C}^m.$$   Let $G_{\epsilon}$ be the isometry group of this form and
write $\mathfrak{g}_{\epsilon}$ for its Lie algebra.   If we set $\mathcal{P}_{\epsilon}(m)$ to be the set of partitions of $m$ in which all even parts occur with even multiplicity when $\epsilon=1$ and all odd parts occur with even multiplicity when $\epsilon = -1$, then
the classification of nilpotent coadjoint orbits takes
the form:

\begin{theorem}[\cite{gerstenhaber}]
The nilpotent coadjoint orbits of $G_{\epsilon}$ are in one-to-one
  correspondence with the set of partitions of $\mathcal{P}_{\epsilon}(m)$.
  \label{theorem:gerstenhaber}
\end{theorem}

When $\epsilon =-1$, $G_\epsilon$ is a group of type $C$; for $\epsilon =1$ it is of type $B$ when $m$ is odd and type $D$ when $m$ is even.  In types $B$ and $C$, the adjoint orbits for the isometry group $G_\epsilon$ coincide with the adjoint orbits for the adjoint group $G_{ad}$ of $\mathfrak{g}_\epsilon$.  In type $D$, however, the adjoint group is $PSO(m/2,\mathbb{C})$ and every $G_\epsilon$-orbit $\mathcal{O}$ that corresponds to a very even partition, that is, a partition with only even parts each of which appears with even multiplicity, is the union of two $PSO(m/2,\mathbb{C})$-orbits. We will denote them as $\mathcal{O}^I$ and $\mathcal{O}^{I\hspace{-.02in}I}$.

In what follows, we will write $\mathcal{O}_\lambda$ for the nilpotent $G_\epsilon$-orbit associated with the partition $\lambda$ when the type of the underlying group is clear.

\subsection{Orbital Varieties in Classical Types}

While partitions were sufficient to describe the set of nilpotent orbits for a classical simple Lie group, somewhat more intricate combinatorial objects are necessary to describe the orbital varieties contained within each nilpotent orbit.

\begin{definition} Let $\lambda$ be a partition.  A {\it Young diagram of shape $\lambda$} is a finite left-justified array of squares the length of whose $i$th row equals the $i$th part of $\lambda$.
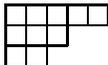
$$
\begin{tiny}
\begin{tableau}
    :.{}.{}.{}.{}.{}\\
    :.{}.{}.{}\\
    :.{}.{}\\
\end{tableau}
\end{tiny}
$$
Write $\mathbb{N}_n=\{1, 2, \ldots, , n\}$.  A {\it standard Young tableau of $\lambda$} is a Young diagram
    of shape $\lambda \dashv n$ whose squares are labeled by elements of $\mathbb{N}_n$ in such a
     way that each element of $\mathbb{N}_n$
      labels exactly one square, and all labels increase along both rows and columns.
\end{definition}

\begin{definition}
Let $r \in \mathbb{N}$ and $\lambda$ be a partition of a positive
integer $m$.
     A {\it standard domino tableau of rank $r$ and shape $\lambda$} is a Young diagram
    of shape $\lambda$ whose squares are labeled by elements of $\mathbb{N}_n \cup \{0\}$ in such a
     way that the integer $0$ labels the square $s_{ij}$ iff $i+j<r+2$, each element of $\mathbb{N}_n$
      labels exactly two adjacent squares, and all labels increase
    weakly along both rows and columns.
We will write $SDT_r(\lambda)$ for the family of all domino tableaux
of rank $r$ and shape $\lambda$  and $SDT_r(n)$ for the family of
all domino tableaux of rank $r$ which contain exactly $n$ dominos.

\end{definition}

Let $\mathfrak{b} \subset \mathfrak{g}$ be a Borel subalgebra,
$\mathfrak{h} \subset \mathfrak{b}$ a Cartan subalgebra, and
$\mathfrak{n}$ the nilradical so that
$\mathfrak{b}=\mathfrak{h}+\mathfrak{n}$.  Recall that for a nilpotent orbit $\mathcal{O}$, the irreducible
components $\Irr(\mathcal{O} \cap \mathfrak{n})$ are
its orbital varieties.  By \cite{joseph}, every orbital variety takes the form $$V(w) = \overline{B(\mathfrak{n}
  \cap w^{-1} \mathfrak{n})} \cap \mathcal{O}$$ for some $w $ in the Weyl group $W$. The set of Weyl group elements which map to the
same orbital variety under this correspondence is known as a {\it
  geometric left cell}. In light of the Robinson-Schensted algorithms which associate elements of $W$ with same-shape pairs of standard Young and domino tableaux (see \cite{garfinkle1}), the following results are somewhat natural:

\begin{theorem}[\cite{joseph}]\label{theorem:ovA}
In type $A$, orbital varieties contained in the nilpotent orbit $\mathcal{O}_\lambda$
are parameterized by the set of standard  Young tableaux of shape
$\lambda.$
\end{theorem}

\begin{theorem}[\cite{mcgovern2},\cite{pietraho:springer}]\label{theorem:ovX}
In types $C$ and $D,$ orbital varieties contained in the nilpotent $G_\epsilon$-orbit
$\mathcal {O}_\lambda$ are parameterized by standard domino tableaux of rank zero and
shape $\lambda$.  In type $B,$ orbital varieties contained in the nilpotent $G_\epsilon$-orbit
$\mathcal {O}_\lambda$ are parameterized by standard domino tableaux of rank one and
shape $\lambda$.
\end{theorem}

The
construction of the Graham-Vogan space associated to an orbital
variety $\calV$ requires us to be able to explicitly identify its
$\tau$-invariant.  We describe how to do this for an orbital variety $\calV_T$
corresponding to a standard tableau $T$.

Let $\Delta$ be the set of roots in $\mathfrak{g}$, $\Delta^+$ the
set of positive roots and $\Pi$ the set of simple roots.  Write
$\mathfrak{g}= \bigoplus_{\alpha \in \Delta^+}
\mathfrak{g}_{-\alpha} \oplus \mathfrak{t} \oplus \bigoplus_{\alpha
\in \Delta^+} \mathfrak{g}_{\alpha}$ for the triangular
decomposition of $\mathfrak{g}$ and let $\mathfrak{b}=\mathfrak{t}
\oplus \bigoplus_{\alpha \in \Delta^+} \mathfrak{g}_{\alpha}$.
Write  $W$ for  the Weyl group, and let $P_\alpha$ be the standard
parabolic subgroup with Lie algebra $\mathfrak{p}_\alpha=
\mathfrak{b} \oplus \mathfrak{g}_{-\alpha}$. Following
\cite{joseph},  for an element $w \in W$, an orbital variety
${\calV}$, and a standard parabolic subgroup we define
\begin{align*}
  \tau (P) & = \{ \alpha \in \Pi \, | \, P_\alpha \subset P\} \text{ and } \\
  \tau (\calV) & = \{ \alpha \in \Pi \, | \, P_\alpha (\calV) = \calV \}.
\end{align*}
We would like to be able to read off $\tau(\calV)$ from the standard
tableau parameterizing $\calV$ as the maximal parabolic subgroup $Q$
stabilizing $\calV$ is precisely the standard parabolic subgroup
satisfying $\tau(Q)=\tau(\calV)$.

\begin{theorem}[\cite{joseph}]
Consider an orbital variety $\calV_T$ in type $A$ that corresponds to the
standard Young tableau $T$ under the above parametrization.  The simple
root $\alpha_i \in \Pi$ lies in $\tau(\calV_T)$ iff  the square labeled
$i$ lies strictly higher in $T$ than the square with label $i+1$.
\end{theorem}

\begin{theorem}[\cite{pietraho:springer}]
Consider an orbital variety $\calV_T$ in type $B,$ $C$, or $D$ that corresponds to the
standard domino tableau $T$ under the above parametrization.  The simple root $\alpha_i$ lies in $\tau(\calV_T)$ iff one
of the following conditions is satisfied:
\begin{enumerate}[\hspace{.2in}(i)]
\item $i=1$ and the domino with label $1$ is vertical,

\item $i>1$ and domino with label $i-1$ lies higher that the domino with label $i$ in $T$.
\end{enumerate}
\label{theorem:tau}
\end{theorem}

Finally, we would like to define a map from the orbital varieties in types $B$, $C,$ and $D$ to orbital varieties in type $A$. Let $\lieg$ be a classical complex Lie algebra of type $X_n=B_n$,
$C_n$, or $D_n$   and let $\mathfrak{n}$ be the unipotent part of
$\mathfrak{b}.$ There is a natural projection map $\pi_A$ from
$\mathfrak{n}$ to $\mathfrak{n}_A,$ the corresponding unipotent part
in type $A_{n-1}$.  Let $\calO$ be a nilpotent orbit of type $X_n$.
The image of an orbital variety for $\calO$  under
$\pi_A$ is always an orbital variety for some nilpotent orbit
$\mathcal{P}$ of type $A$. In fact, if $\mathcal{P}$ arises in this
way, then {\it all} of its orbital varieties lie in the image of
$\pi_A$ for $\calO.$  To describe this in terms of the underlying combinatorics, we need a result of Carre and Leclerc.

\begin{theorem}[\cite{cl}]
There is a bijection
$$
\begin{CD}
SDT(\lambda) @>>{(\pi_1,\pi_2)}>  \amalg_{\nu} Yam_2(\lambda, \nu)
\times SYT(\nu). \end{CD}$$
where $ Yam_2(\lambda, \nu)$ is the family of Yamanouchi domino tableaux of shape $\lambda$ and evaluation $\mu$.
\end{theorem}
The bijection itself is an algorithm that takes a tableau $T$ and modifies
it successively until its column reading becomes a Yamanouchi word.
The standard Young tableau  records the
sequence of moves.  We are interested only in the second coordinate of this map.

\begin{definition}
Define a map $$\pi_A: SDT(n) \longrightarrow SYT(n)$$ by $\pi_A (T)
= \pi_2 (T)$ where $\pi_2$ is the second component of the
Carre-Leclerc map.  We also denote by $\pi_A$ the map induced on
orbital varieties  obtained by identifying $T$ with $\calV_T$.
\end{definition}


\section{Restriction to Spherical Orbital Varieties}\label{section:restrictiontospherical}

Armed with a description of the orbital varieties contained in a
given nilpotent orbit as well as the corresponding
$\tau$-invariants, we now begin to describe the Graham-Vogan
representations attached to a nilpotent orbit in the setting of
classical groups.  We begin by illustrating our method with an  example, which
is sufficiently na\"{\i}ve to quickly describe our approach.

\subsection{Model Example}

We will calculate the infinitesimal character associated to
$V(\calV, \pi)$
 constructed from a particular orbital variety in type $C$.
Suppose $G=Sp(8)$ and realize the Lie algebra $\mathfrak{g}$
as a set of $8 \times 8$ matrices of the form
$$\mathfrak{sp}(8)= \bigg\{\mathfrak{m}(A,B,C)=
\begin{tiny}\left(
\begin{array}{cc}
A & B \\
C & -A^t
\end{array}
\right)\end{tiny}| \, A, B, C \in M_4(\mathbb{C}) \text{ and } B, C \in
Sym_4(\mathbb{C}) \bigg\}.$$

Let $\calO$ be the nilpotent coadjoint
orbit in $\mathfrak{g}^*$ corresponding to the partition
$[2^3,1^2].$ It has dimension 18. There are four orbital varieties
contained in $\calO,$ corresponding to the domino tableaux:
\begin{tiny}
$$\begin{tableau}
:>1\\
:^2^3\\
:;\\
:^4\\
\end{tableau}
\hspace{0.99in}
\begin{tableau}
:^1^2\\
:;\\
:>3\\
:^4\\
\end{tableau}
\hspace{0.99in}
\begin{tableau}
:>1\\
:^2^4\\
:;\\
:^3\\
\end{tableau}
\hspace{0.99in}
\begin{tableau}
:>1\\
:>2\\
:>3\\
:^4\\
\end{tableau}.
$$
\end{tiny}

Let $\calV$ be the orbital variety corresponding to
the first domino tableau. Then $\dim \calV = \smash{\frac{1}{2}}
\dim \calO = 9.$ As a representative, we take
$f=\mathfrak{m}(A,B,0)$ with
$$A=
\begin{tiny}
\left(
\begin{array}{cccc}
0 & 0 & 0 & 0 \\
0 & 0 & 1 & 0\\
0 & 0 & 0 & 0 \\
0 & 0 & 0 & 0
\end{array}
\right)
\end{tiny} \text{\hspace{.2in} and \hspace{.2in}} B= \begin{tiny}\left(
\begin{array}{cccc}
0 & 0 & 0 & 0 \\
0 & 0 & 0 & 0 \\
0 & 0 & 0 & 0 \\
0 & 0 & 0 & 1
\end{array}
\right)
\end{tiny}.$$ To describe the
Graham-Vogan space for $\calV$, we need to compute the following parameters:
an admissible orbit datum $(\pi, \mathcal{H}_{\pi}),$  $\calV^o,$ the smooth part of $\calV,$
 the stabilizing parabolic $Q_\calV \subset G,$
a smooth representation $(\gamma, W_{\gamma})$ of $Q_{\calV},$
and a $G$-equivariant isomorphism of vector bundles $j_{\pi},$
where notation  is as in Section \ref{section:construction}

Write $G_f$ for the isotropy subgroup of $f$ and $ \mathfrak{g}_f$
for its Lie algebra. As $G$ is complex, the metaplectic cover
$\widetilde{G}_f$ is isomorphic to $G_f  \times
\mathbb{Z}/2\mathbb{Z}$.  We choose one admissible orbit datum; it
is trivial on $G_f^\circ$ and acts by the non-trivial character on
$\mathbb{Z}/2\mathbb{Z}$.
The orbital variety $\calV$ is smooth so that in the notation of the
first section, $\calV^o=\calV$.  From Theorem \ref{theorem:tau}, we
find that the stabilizer of $\calV$ is the standard parabolic
subgroup $Q_{\calV}$ with Levi factor isomorphic to $GL(2) \times
GL(2)$.  One can quickly check that, in this case, both $Q_{\calV}$
and the standard Borel subgroup $B$ act with dense orbit on $\calV$.

This observation simplifies calculations, as it allows us to replace
the Lagrangian covering $G \times_{Q_{\calV}} \calV$ by  $G/Q_f,$
where $Q_f = Q_{\calV} \cap G_f$ and $ \calV$ contains
$Q_{\calV}/Q_f$ as a dense subset.  We note that  $B/B_f$ is also
dense in $\calV$. The equivariant line bundle $\tau^* \calV_{\pi}$
is induced by a character $\alpha$ of $B_f$.  It is given by the
square root of the absolute value of the real determinant of $B_f$
acting on the tangent space $\mathfrak{b}/\mathfrak{b_f}$ of $\calV$
at $f$.  This is
$$\alpha
\left(
\begin{array}{cc}
A & * \\
0 & A^{t^{-1}}
\end{array}
\right) = | \, t_1^3 \, t_3^6 \,|^{-1} \text{ , where } A =  \begin{tiny}\left(
\begin{array}{cccc}
t_1 & * & * & * \\
0 & t_3 & * & * \\
0 & 0 & t_3 & * \\
0 & 0 & 0 & 1
\end{array}
\right)\end{tiny}.$$ Because we are looking for a map $j_{\pi}$, we would like
to find a character $\gamma$ of $B$ whose restriction to $B_f$ is
$\alpha.$  Such a character is given by
$$ \gamma
\left(
\begin{array}{cc}
A & * \\
0 & A^{t^{-1}}
\end{array}
\right) = | \, t_1 \, t_2  \, t_3 \, t_4 \,|^{-3} \text{ , where } A
= \begin{tiny}\left(
\begin{array}{cccc}
t_1 & * & * & * \\
0 & t_2 & * & * \\
0 & 0 & t_3 & * \\
0 & 0 & 0 & t_4
\end{array}
\right)\end{tiny}.$$ The character $\gamma$ extends uniquely from $B$ to
$Q_{\calV}$. Let the half-density bundle on $G/Q_{\calV}$ be given
by the character $\rho_{Q_{\calV}}$ and define another character
$\gamma'$ on $Q_{\calV}$ to equal $\gamma \otimes
\rho^{-1}_{Q_{\calV}}.$  Then
$$V(\calO, \calV, \pi, \gamma, j_{\gamma, \pi}) \subset
  Ind_Q^G(\gamma').$$
Hence the infinitesimal character that we associated to the
representation space $V(\calV, \pi)$ equals
$-\smash{(\frac{3}{2},\frac{3}{2},\frac{3}{2}, \frac{3}{2}) + \rho =
(\frac{5}{2},\frac{3}{2},\frac{1}{2}, \frac{1}{2})},$ where
$\rho=(4,3,2,1)$ and equality is up to Weyl group action.  This is
precisely the unique infinitesimal character attached to the orbit
$\calO_{[2^3,1^2]}$ by McGovern in \cite{mcgovern}.  Similar calculations for
the other orbital varieties in this orbit yield the same
infinitesimal character.

A significant simplification in this example came from the fact that
the parabolic subgroup $Q_{\calV}$ acted with dense orbit on
$\calV.$ It made it easy to find the isomorphism $j_{\pi}.$
Unfortunately, this is not always the case.

\begin{example}[\cite{melnikov}]
Let $G = SL_9$ and let
$$\raisebox{2ex}{T= \;}
\begin{tiny}
\begin{tableau}
:.1.2.3.6.9\\
:.4.5.8\\
:.7\\
\end{tableau}
\end{tiny}
$$
Then $\calV_T,$ the
  orbital variety in $\calO_{[5,3,1]}$ corresponding to $T$ has
  dimension $31.$  However, $\dim Q_{\calV} \cdot f \leq 30$ for all
  $f \in \calV.$
\end{example}

This example can be extended to produce other instances where $Q_{\calV}$ does not act with dense orbit on
$\calV$ both, in larger groups of type $A$ as well in groups of other classical types.  We will restrict our attention to a class of nilpotent orbits all of whose orbital varieties do admit a dense orbit of their stabilizing parabolic, but note that every nilpotent orbit contains at least one orbital variety with this property.

\subsection{Spherical Orbital Varieties and Orbits of $S$-type}

We would like to use the methods of our model example to calculate
the infinitesimal character associated to $V(\calV, \pi)$ for as
many nilpotent orbits as feasible.  The main assumption required is
that the stabilizer of an orbital variety has a dense orbit in that
variety. Such orbital varieties are called of $S${\it -type}, as are
the nilpotent orbits {\it all} of whose orbital varieties satisfy
this condition. Among classical groups, there is a class of {\it
small} nilpotent orbits that are of $S$-type.  We first describe
this set and then place it among other  important nilpotent
coadjoint orbits.

Let $G$ be a complex simple Lie
group and $B$ a Borel subgroup.  We will say that
a nilpotent coadjoint orbit $\calO \subset \mathfrak{g}^*$ is
{\it spherical} iff it contains an open $B$-orbit.
The work of D.~Panyushev provides a concise description of spherical
nilpotent orbits contained in classical groups.

\begin{theorem}[\cite{panyushev}]
A nilpotent orbit $\mathcal{O}_\lambda$  of a complex classical simple Lie group that is parametrized by the partition $\lambda$ is spherical if $\lambda$ is of the form:
\begin{enumerate}[\hspace{.2in}(i)]
  \item $[2^b, 1^c]$ in type $A,$

  \item $[3^a,2^{2b},1^c]$ with $a \leq 1$ in type $B,$

  \item $[2^b, 1^{2c}]$ in type $C,$ and

  \item $[3^a,2^{2b},1^c]$ with $a \leq 1$ in type $D.$
\end{enumerate}
\end{theorem}

A few properties characterize spherical orbits.   First, they are
precisely the orbits which contain a representative that is a sum a
root vectors corresponding to orthogonal simple roots (\cite{mcgovern3} and \cite{panyushev2}).
Furthermore, for complex simply-connected semisimple Lie groups, there is an orbit for which the $G$-module structure of $R(\calO)$, the coordinate ring of regular functions, has all multiplicities either $0$ or $1$ \cite{mcgovern3}. The largest such orbit is called the {\it model
orbit}.  Spherical orbits may be characterized  as those nilpotent
coadjoint orbits contained in the closure of the model orbit.

\begin{theorem}[\cite{mcgovern3}]
Let $\epsilon = 0$ or $1$.
In each of the classical types, the model orbit is the largest
spherical nilpotent orbit and is parametrized by the following
partition:
\begin{enumerate}[\hspace{.2in}(i)]
  \item $[2^n, 1^{\epsilon}]$ in type $A_{2n+\epsilon-1}$,

  \item $[3,2^{4m-2\epsilon},1^{2\epsilon}]$ in type $B_{2(2m-\epsilon)+1},$

  \item $[2^n]$ in type $C_{2n},$ and

  \item $[3,2^{2m-2},1^{1+2\epsilon}]$ in type $D_{2(2m+\epsilon)}.$
\end{enumerate}
\end{theorem}

Following A.~Melnikov in \cite{melnikov}, we will say that an orbital variety $\calV \subset \calO$ is of {\it $S$-type} iff it
admits a dense orbit of its maximal stabilizing parabolic  $Q_{\calV}$ and extend the terminology to nilpotent coadjoint orbits all of whose
 orbital varieties are of $S$-type.  The following proposition is a consequence of the dimension argument in Corollary
\ref{cor:dimension}.

\begin{proposition}
In the setting of complex classical simple Lie groups, all spherical nilpotent orbits are of $S$-type.
\end{proposition}

Although we will restrict our attention to
spherical nilpotent orbits,  for
completeness, we provide a partial description of the $S$-type orbits in groups
of type $A$.  The result is incomplete, as it fails to resolve the status of a
number of nilpotent orbits.

\begin{theorem}[\cite{melnikov}]
A nilpotent orbit $\calO_{\lambda}$ in type $A$ is of $S$-type
whenever $\lambda$ satisfies one of the following:
\begin{enumerate}[\hspace{.2in}(i)]
  \item $\lambda > (n-4,4),$
  \item $\lambda = (\lambda_1, \lambda_2, 1, \ldots, 1)$ with
    $\lambda_2 \leq 2,$
  \item $\lambda = (2, \ldots)$ where $\lambda_i \leq 2$ for all $i.$
\end{enumerate}
If we suppose that $n \geq 13$, the partition $\lambda$ has $\lambda_2
>2$, and
$(5,3,1, \ldots) \leq \lambda \leq (n-4,4)$
in the usual partial order on partitions,  then the orbit
$\calO_{\lambda}$ in type $A_{n-1}$ is not of $S$-type.
\end{theorem}

We finish this section by listing how
spherical orbits fit among two other important classes of nilpotent
orbits.  A nilpotent orbit is {\it rigid} if it is not
induced from any proper parabolic subalgebra.  It is {\it special} if
it is in the range of a particular order-reversing map $d$, see for instance \cite{cm}(6.3.7)
for a characterization.  The following propositions are immediate consequences of the results in \cite{cm}.

\begin{proposition}
All nilpotent orbits are special in type $A.$  In the other classical types, the nilpotent orbit
$\mathcal{O}_\lambda$   that is parametrized by the partition $\lambda$ is spherical and special if $\lambda$ is of the form:
\begin{enumerate}[\hspace{.2in}(i)]
  \item
    $[3,2^{2b},1^c]$ or $[1^c]$ in type $B,$
  \item $[2^{2b},1^{2c}]$ or $[2^b]$ in type $C,$
  \item $[2^{2b},1^c]$ or$[3,1^c]$ in type
    $D.$
\end{enumerate}
\end{proposition}

\begin{proposition} All non-zero orbits are not rigid in type $A.$ In the other classical types, the nilpotent orbit
$\mathcal{O}_\lambda$   that is parametrized by the partition $\lambda$ is spherical and non-rigid if $\lambda$ is of the form:
\begin{enumerate}[\hspace{.2in}(i)]

  \item $[3,1^{2c}]$ or $[2^{2b},1^2]$ in type $B,$

  \item  $[2^2,
    1^{2c}]$ or $[2^{2c}]$  in type $C,$ and

  \item
    $[3,1^{c}]$ or $[2^{2c}]$ in type $D.$
\end{enumerate}
\end{proposition}


\subsection{Basepoints in Orbital Varieties}\label{subsection:basepoint}

From the previous section, we know that each spherical orbital
variety $\calV$ contains a point whose orbit under the Borel
subgroup is dense in $\calV$. We would like a
simple expression for some such point to simplify the forthcoming
calculations. For orbital varieties within classical nilpotent
orbits, such an expression can be easily read off from the standard
tableau corresponding to $\calV$.

In type $A$, such a basepoint is essentially defined in
\cite{melnikov}. We provide a slightly more general construction and
extend the result to other classical types.  The main tool for the
latter is the surjection $\pi_A$ from domino tableaux onto standard Young
tableaux defined in \cite{cl}.  It induces a map on the level of
orbital varieties that helps us define the basepoint in the ``type
$A$ component'' of each orbital variety.

\subsubsection{Type A}

Consider a spherical nilpotent orbit $\calO$ and let $\calV_T
\subset \calO$ be the orbital variety associated to the standard
Young tableau $T \in YT(n)$. Let $T^i$ denote the set of labels
contained in the $i$-th column of $T,$  so that in our case $T^i =
\emptyset$ if $i > 2$.  We will define a point $f_T$ contained in
$\calV_T$ whose orbit under the Borel subgroup is dense in
$\calV_T.$   For $x\in \mathbb{N}$, let $\tilde{x} = n+1 -x$ and let us adopt the notation from \cite[IV.1]{knapp:example}, writing $E_{e_i-e_j}$ for the matrix with the $ij$-entry equal to one and zero otherwise.

\begin{proposition}\label{proposition:basepointA}
Let $\phi: T^2 \longrightarrow T^1$ be an injection with the
property that $\phi(k) < k $ for all $k \in T^2$. Such a map always
exists, and furthermore, the point
$$f_T = \sum_{k \in T^2} E_{e_{\tilde{k}} - e_{\widetilde{\phi(k)}}}$$
is contained in the variety $\calV_T$.
\end{proposition}

\begin{proof}
The fact that a map $\phi$ always exists is clear by inspection.  A
spherical nilpotent orbit in type $A$ is uniquely determined by the
rank of its elements.  For each $f_T$ defined above, $f_T^2 =0$, so
it lies in {\it some} spherical orbital variety.  That it lies
precisely in $\calV_T$ follows from induction and the above rank
condition.
\end{proof}

This definition includes Melnikov's construction as a special case.
More precisely, it is always possible to  choose $\phi$ in such a
way so that $\phi(k)=k-1$ whenever $\alpha_{\widetilde{{k-1}}}
\notin \tau(T)$.  In this incarnation, $f_T$ is a {\it minimal
representative} of $\calV_T$ in the sense described below. Let $f
\in \mathfrak{n}$ and for its root space decomposition, let us write
$$f= \sum_{\epsilon \in \Delta^+} c_{\epsilon}(f) E_{\epsilon}.$$

\begin{definition}
An element $f \in \calV$ is a {\it representative} of $\calV$ if $f$ does
not belong to any other orbital varieties. A representative $f$ of
$\calV$ is {\it minimal} if
\begin{enumerate}[\hspace{.2in}(i)]
  \item each $c_{\epsilon}(f) \in \mathbb{Z}$,
    \item for every $\alpha_i \notin \tau(\calV)$, $c_{\alpha_i}(f) \neq 0$,
   \item If $g \in \calV$ also satisfies the
        above, the the number of non-zero $c_{\epsilon}(g)$ will be greater than or equal to
        the number of
        non-zero $c_{\epsilon}(f)$.
 \end{enumerate}
\end{definition}

We would like the basepoints we choose to be minimal
representatives.  In type $A$, we have already
seen that this is always possible and in further work we would like
$f_T$ to be close to satisfying this condition.

\begin{example}
Consider the orbital variety $\calV_T$ associated with
  the standard Young tableau
$$\raisebox{3ex}{T = \;}
\begin{tiny}
\begin{tableau}
:.1.3\\
:.2\\
:.4\\
:.5\\
\end{tableau}
\end{tiny}
$$
The points $f_1 = E_{e_1-e_2}$ and $f_2=E_{e_1-e_3}$
both lie in $\calV_T$ and both are $f_T$ for some choice of injection $\phi$ by Proposition \ref{proposition:basepointA}. However, only $f_1$ is a minimal
representative of $\calV_T$.

\end{example}

\subsubsection{Other Classical Types}

Let $X =B$, $C$, or $D$, and let $\calV_T$ be the orbital variety in a
spherical nilpotent orbit of type $X$ associated with the standard domino
tableau $T.$  In search of a suitable basepoint, we first define a matrix $M_T^X$ from the horizontal dominos of $T$.  Let $N^T$ to be the set of labels of the horizontal dominos in $T$ and define
$S^T$ be the subset of $N^T$ whose underlying dominos intersect the first column
of $T$.
If $M$ is a family of sets of integers, let $M^\circ$ denote the
union of all integers contained in elements of $M.$   Write $T(m)$ for the domino tableau consisting of the first $m$ dominos of $T$ and $D(m)$ for the domino with label $m$. We now
inductively define a set $N_1^T$ of {\it pairs} of labels in $T$
by $N_1^\varnothing = \varnothing$ and
$$
N_1^T = \left\{ \begin{array}{ll}
    N_1^{T(n-1)} \cup \{\{k,n\}\} & \text{ if } D(n) \in S^T \setminus
    (N_1^{T(n-1)})^\circ \\
                                  & \; \text{ and if } X=C,\ k=n-1, \\
    N_1^{T(n-1)} & \text{ otherwise.} \\
\end{array} \right.
$$
Finally, let $ N_2^T = S^T \setminus (N_1^T)^\circ \text{ and }
N_3^T = N^T \setminus ((N_1^T)^\circ \cup N_2^T).$  Note that $N_3^T$ is always empty in type $C$ while $N_2^T$ is always
empty in types $B$ and $D.$

\V
\begin{example}
Suppose $T$ and $U$ are the following domino tableaux and note that $U$ can be viewed as a domino tableau of type $C$ as well as $D$:
$$
\raisebox{3ex}{$T=$ \;}
\begin{tiny}
\begin{tableau}
:.0>1\\
:>2\\
:^3^4\\
:;\\
:>5\\
\end{tableau}
\end{tiny}
\hspace{1in} \raisebox{3ex}{$U=$ \;}
\begin{tiny}
\begin{tableau}
:>1\\
:^2^3\\
:;\\
:>4\\
\end{tableau}
\end{tiny}
\raisebox{3ex}{   .}
$$
For the tableau $T$ of type $B$,
$S^T = \{2,5\} \text{, } N_1^T = \{\{2,5\}\} \text{, } N_2^T = \varnothing \text{, and } N_3^T = \{1\}.$
For the tableau $U$,
$S^U =  \{1,4\} \text{, }  N_1^U = \varnothing \text{, }  N_2^U =  \{1,4\} \text{, and } N_3^U = \varnothing$ when it is viewed as a domino tableau of type $C$, and
$S^U = \{1,4\} \text{, }N_1^U =  \{\{1,4\}\} \text{, } N_2^U =  \varnothing  \text{, and } N_3^U = \varnothing$ when it is viewed as a domino tableau of type $D.$
\end{example}


As in the previous section,  we will adopt the notation for simple roots from \cite[IV.1]{knapp:example}, write $E_\alpha$ for a basis vector for the root space $\mathfrak{g}_\alpha$ of a simple root $\alpha$ and take $T_i \in \mathfrak{t}_i$ where $\mathfrak{t} = \oplus_{i \leq rank \, G }\mathfrak{t}_i.$ Let $$M_T^X
= \sum_{\{i,j\} \in N_1^T} E_{e_{\tilde{i}} + e_{\tilde{j}}} +
\sum_{\alpha \in U_T^X } E_{\alpha}, $$
where $U_T^X$ is the set of roots defined by
$$
U_T^X = \left\{ \begin{array}{ll}
\{e_{\widetilde{k-1}}+e_{\tilde{k}} \, | \, k \in N_3^T \} &\text{if } X=D, \\
\{2e_{\tilde{k}} \, | \, N_2^T \} & \text{if } X=C, \text{ and}  \\
\{e_{\tilde{k}} \, | \, k \in N_3^T \} \cup \{e_{\tilde{3}} \, | \,
2
\in T^3 \text{ and } 3 \in T^2 \} & \text{if }X=B. \\
\end{array}
\right.
$$

\begin{definition}
For $X= B, C,$ or $D$, and a domino tableau $T$, let
$$f_T^X = f_{\pi_A(T)} + M_T^X$$
where $f_{\pi_A(T)}$ is a minimal representative of the type $A$ orbital variety $\mathcal{V}_{\pi_A(T)}$ interpreted naturally as lying inside the Lie algebra
of type $X$. \label{definition:basepoint}
\end{definition}

\begin{proposition}\label{proposition:basepointX}
The point $f_T^X$ is a minimal representative of $\calV_T.$
\end{proposition}

\begin{proof}
The proof is a little simpler if we use an alternate parametrization of orbital varieties which uses the set of admissible domino tableaux with signed closed clusters
$\Sigma DT_{cl}$.  We refer the reader to \cite{pietraho:springer} for the details, and let $\Phi$ be the bijection between the two parameter sets defined therein.

Define $T' = \Phi^{-1}(T) \in \Sigma DT_{cl}(shape \, T)$. We first
show that  $f_T^X \in \calV_S$, where $S= \Phi(T^*)$ and $T^* \in
\Sigma DT_{cl}(shape \, T)$ has the same underlying domino tableau
as $T'$.  We then show that $T'$ and $T^*$ must share the set of
closed clusters with positive sign, which implies that $S=T$ by the
definition of $\Phi$.  This verifies that $f_T^X$ is a representative
of $\calV_T$.  Minimality of $f_T^X$ may then be checked by
inspection.

We would like to show that for all $k \leq n,$ $f_{T(k)}^X \in
\calO_{\text{shape } T'(k)}.$  By induction, it is enough to verify
this for $k = n-1.$  Note that for spherical orbits, the partition
of the orbit containing a nilpotent element $f$ is completely
determined by rank $f$ and rank $f^2$.  The above statement can be
now verified  by inspecting the definition of $f_T^X$ and comparing
rank $f_{T(n-1)}^X$ and rank $(f^X_{T(n-1)})^2$ with rank $f_T^X$ and rank
$(f_T^X)^2$. In this way, we have  $f_T^X \in \calV_S$, where $S=
\Phi(T^*)$ and $T^*$ is some tableau in $\Sigma DT_{cl}(shape \, T)$
sharing its underlying tableau with $T'.$

Now note that if $\calC$ is a closed cluster of $T'$ or $T^*$, then
because the orbit $\calO_{shape \, T}$ is spherical, the initial
cycle $I_\calC$ through $\calC$ must have the form $I_\calC = \{i,
i+1, \ldots , j\}$. Theorem \ref{theorem:tau} implies that the
simple root $\alpha_{\tilde{i}} \in \tau(T)$ iff there is a closed
cluster $C \in \calC^+$ with $I_\calC = \{i, i+1, \ldots , j\}$ for
some $j$.  Further note that if $C \in \calC^+$, then  $E_{e_i +
e_j}$ appears in the expansion of $f_T^X$ with non-zero coefficient
while $E_{e_i - e_j}$ has coefficient zero.  Similarly, if $C \in
\calC^-$, then $E_{e_i - e_j}$ appears in the expansion of $f_T^X$
with non-zero coefficient while $E_{e_i - e_j}$ has coefficient
zero. But this forces $\Phi(T^*)$ to have the same $\tau$-invariant
as $\Phi(T')$, which implies that $\Phi(T^*)=\Phi(T').$  Hence $f_T^X$
is a representative of $\calV_T$.
\end{proof}

\begin{lemma}
  Consider an orbital variety $\calV_T$ in a spherical nilpotent orbit
  of classical type that corresponds to the standard tableau
  $T$, and let $ Q= Q_{\calV_T} \supset B$ be the maximal parabolic stabilizing
  it.  Then the orbits
  $B \cdot f_T^X \text{ and } Q \cdot f_T^X $
  are dense in $\calV_T.$
\label{lemma:dense}
\end{lemma}

\begin{proof}
  For the result in  in type $A$, see \cite[4.13]{melnikov}.
  The present result follows by induction from Corollary \ref{cor:dimension} below.
\end{proof}

\begin{example}
  Let $X=C$ and consider the orbital variety $\calV_T$ associated with
  the domino tableau
$$
\raisebox{3ex}{T= \; }
\begin{tiny}
\begin{tableau}
:>1\\
:^2^3\\
:;\\
:>4\\
:>5\\
\end{tableau}
\end{tiny}
\raisebox{3ex}{\text{ .  Then the Young tableau  $\pi_A(T)= \; $}}
\begin{tiny}
\begin{tableau}
:.1.3\\
:.2\\
:.4\\
:.5\\
\end{tableau}
\end{tiny}
\raisebox{3ex}{   .}
$$
We have  $N =\{1,4,5\}$, $N_1
  =\{1\}$ and $N_2 = N \setminus N_2$.
Finally, the basepoint $f_T^C = \begin{tiny}\left(
\begin{array}{cc}
A & M \\
0 & -A^t \\
\end{array}
\right)\end{tiny},$ where $$A= f_{\pi_A(T)} =
\begin{tiny}
\left(
\begin{array}{ccccc}
0 & 0 & 0 & 0 & 0 \\
0 & 0 & 0 & 0 & 0 \\
0 & 0 & 0 & 1 & 0 \\
0 & 0 & 0 & 0 & 0 \\
0 & 0 & 0 & 0 & 0 \\
\end{array}
\right)
\end{tiny}
 \text{ \hspace{.1in} and \hspace{.1in}} M=
 \begin{tiny}\left(
\begin{array}{ccccc}
0 & 1 & 0 & 0 & 0 \\
1 & 0 & 0 & 0 & 0 \\
0 & 0 & 0 & 0 & 0 \\
0 & 0 & 0 & 0 & 0 \\
0 & 0 & 0 & 0 & 1 \\
\end{array}
\right)
\end{tiny}.$$
\end{example}


\subsection{Induction}

Our calculation of infinitesimal characters of Graham-Vogan
representations attached to the orbital variety $\calV_T$ will
proceed by a type of induction on the standard tableau $T.$  As in
our model example, we would like to describe the action of $Q_f$ on
the space $\mathfrak{q}/\mathfrak{q_f}.$  Because we consider only
spherical orbits, it necessary only to describe this action on $\mathfrak{b}/\mathfrak{b_f}.$  In this section, we
describe it, verifying Lemma \ref{lemma:dense}
in the process.

Fix a standard tableau $T$ of a given classical type and write
$\calV_T$ for the orbital variety corresponding to it.  Ideally, we
would like to be able extract information about
$\calV_T$ from $\calV_{T(n-1)}$ and in this
manner set up a type of induction.  However, the standard domino
tableau $T(n-1)$ does not always correspond to an orbital variety of
the same classical type as $\calV_T$, so in order for induction to
make sense, we have to be careful.   To this effect, we define a  standard tableau
$T^{\downarrow}$ by
$$
T^{\downarrow} = \left\{
\begin{array}{ll}

MT(\mathcal{C}, T(n-1)) & \text{$X=B$ or $D$, type $\calV_{T(n-1)} \neq X$ and}\\

& \text{ $\mathcal{C}$ the cycle in $T(n-1)$ through $n-1$,}\\

T(n-2) & \text{$X=C$, $D(n)$ and $D(n-1)$ are horizontal,} \\

& \text{ while $D(n-2)$ is not,}\\
T(n-1) & \text{otherwise.}
\end{array}\right.
$$
The notion of a cycle and the moving-through map $MT$ are defined in \cite[\S 5]{garfinkle1}.
With this definition, $\text{shape } T ^\downarrow$ and
$\text{shape } T$ are partitions of the same classical type.
Therefore, we are able to associate an orbital variety
$\calV_{T^\downarrow}$ of the same type as $\calV_T$ to the standard
tableau $T^\downarrow$. We will write $f^\downarrow$ for
$f_{T^\downarrow}^X$, and $\mathfrak{b}^\downarrow,
\mathfrak{q}^\downarrow$, and $\mathfrak{g}^\downarrow$ for the  Lie
algebras corresponding to $\calV_{T^\downarrow}$.

As in our model example, we would like to describe the action of
$Q_f$ on the space $U= \mathfrak{b}/\mathfrak{b}_{f_T}$.  If we
think inductively, however, we can break this task down into a study
of the quotients
$U_{n} = (\mathfrak{b}/\mathfrak{b}_{f_T}) /
(\mathfrak{b}^\downarrow /
\mathfrak{b}^\downarrow_{f^\downarrow}).$
It will be often convenient to divide our work into cases that arise
from an inductive construction of of the representative $f_T^X$. Let $\iota$ be the natural inclusion map of the Lie algebra of type $X$ of rank $n-1$ to the one of rank $n$.  The
cases are distinguished by the possible forms of the difference $f_T^X
- \iota(f_{T^\downarrow}^X).$ We describe the possibilities along with
what they imply on the level of tableaux.
\begin{itemize}
\item[(C1)] When $f_T^X = \iota(f_{T^\downarrow}^X),$  the
domino  $T \setminus T^\downarrow$ lies entirely in the first column of $T$.
\item[(C2)]  When $f_T^X = \iota(f_{T^\downarrow}^X) + E_{e_1-e_{\tilde{\phi}(n)}}$, the
 domino $T \setminus T^\downarrow$ lies entirely in the second column
  of $T$.
\item[(N1)]  When $f_T^X = \iota(f_{T^\downarrow}^X)  + E_{e_1+e_{\tilde{k}}}$ and $X=B$ or
 $D$,  then this is the case when $T^\downarrow \neq T(n-1)$ and $\{k, k+1, \ldots, n-1\}$
 is a cycle in $T(n-1)$.  If $X=C$ and $\tilde{k} \neq 2$, then $k=n-1$ and
  $T^\downarrow=T(n-2)$.

\item[(N2)] When $f_T^X = \iota(f_{T^\downarrow}^X) + E_{2e_1} $, then
  $X=C$ and $T^\downarrow = T(n-1)$.

\item[(N3)] When $f_T^X = \iota(f_{T^\downarrow}^X) + E_{e_1-e_{2}} + E_{e_1},$ we
have $X=B$ and $T \setminus T^\downarrow$ is a horizontal domino that intersects
the third column of $T$. When $f_T = \iota(f_{T^\downarrow}) + E_{e_1-e_{2}} + E_{e_1+e_2},$
we have $X=D$ and again $T \setminus T^\downarrow$ is a horizontal domino that intersects the
 third column of $T$.

\item[(*)] When $f_T^X = \iota(f_{T^\downarrow}^X)+ E_{e_1} $, we have $X=B$ and
$T \setminus T^\downarrow = D(3) \in T^2$ while $D(2) \in N_3^T$.
\end{itemize}

\begin{lemma}
Consider a standard tableau $T$, the orbital variety $\calV_T$, and representative
$f_T$ or $f_T^X.$ Let $\phi$ be the injection from Propositions \ref{proposition:basepointA} and \ref{proposition:basepointX} used to define this representative.  In each of the above cases, the space $U_n$ is:

\begin{itemize}
\item[(C1)] $U_n = \bigoplus_{T^2}
      \mathfrak{g}_{e_1-e_{\tilde{i}}}$  in type $A$.  In the other classical types, let $N= (N_1^T)^\circ \cup N_3^T \cup \{3\}$ in type $B$,  $N= (N_1^T)^\circ \cup N_3^T$ in type $D$ and $N= (N_1^T)^\circ \cup N_2^T$ in type $C$.  Then

$$U_n= \bigoplus_{N}
      \mathfrak{g}_{e_1-e_{\tilde{i}}} \oplus \bigoplus_{(\pi_A (T))^2}
      (\mathfrak{g}_{e_1-e_{\tilde{i}}} \oplus \mathfrak{g}_{e_1+e_{\widetilde{\phi(i)}}})$$

\item[(C2)]
Define $V$ and $W$ by
\begin{align*}
       V  =  &  \bigoplus_{\substack{j>\widetilde{\phi(n)}\\ \tilde{j}
          \notin N^T \cup (\pi_A(T))^2}}
        \mathfrak{g}_{e_{\tilde{\phi}(n)}-e_j} \oplus
        \bigoplus_{\substack{(\pi_A (T(n-1)))^2
                   \\ \phi(i) > \phi(n)}}
        \mathfrak{g}_{e_1-e_{\tilde{i}}} \oplus \mathfrak{t}_1
       \\
W = &    \bigoplus_{\substack{j \neq \tilde{\phi}(n) \cup 1 \\ j
\neq
            \tilde{\phi}(i), \tilde{i} \in (\pi_A (T^\downarrow))^2}}
       \mathfrak{g}_{e_{\tilde{\phi}(n)}+e_j}. \\
\end{align*}
Then $U_n =V$ in type $A$ and $U_n=V \oplus W \oplus \mathfrak{g}_{e_1}
      \oplus_{N_3^T = \varnothing} \mathfrak{g}_{e_{\tilde{\phi}(n)}}$ in type $B$.  In type $C$,
$$U_n=V \oplus W \oplus
        \mathfrak{g}_{2e_{\tilde{\phi}(n)}} \oplus
        \mathfrak{g}_{e_1+e_{\tilde{\phi}(n)}}, $$ while in type $D$,  $$U_n=V \oplus W \oplus_{N_3^T}
      (\mathfrak{g}_{e_1-e_{\tilde{i}}} \oplus \mathfrak{g}_{e_1 -
        e_{\widetilde{i-1}}}).$$

\item[(N1)] Let  $N_1^T=\{\{k,n\}\}$  and write $N=(N_1^T)^\circ$.  Then
\begin{equation*}
U_n=  \bigoplus_{N} \mathfrak{g}_{e_1-e_{\tilde{j}}}
      \oplus \bigoplus_{(\pi_A (T^\downarrow))^2}
      ( \mathfrak{g}_{e_1-e_{\tilde{i}}} \oplus \mathfrak{g}_{e_1+e_{\tilde{\phi}(i)}}) \oplus \mathfrak{t}_1
\end{equation*}

\item[(N2)] This case arises only in type $C$.  Write $N= N^{T^\downarrow} \setminus \{1\}$.  Then
\begin{align*}
U_n = & \phantom{W}
  \bigoplus_{N}(\mathfrak{g}_{e_1-e_{\tilde{i}}} \oplus
    \mathfrak{g}_{e_2-e_{\tilde{i}}}) \oplus\\
 &  \bigoplus_{\substack{(\pi_A
      (T^\downarrow))^2}} (\mathfrak{g}_{e_1-e_{\tilde{i}}} \oplus
  \mathfrak{g}_{e_2-e_{\tilde{i}}} \oplus \mathfrak{g}_{e_1+\tilde{\phi}(i)} \oplus
  \mathfrak{g}_{e_2 +\tilde{\phi}(i)}) \oplus \mathfrak{t}_1.
\end{align*}

\item[(N3)] This case arises only in types $B$ and $D$.  In the former case, $$U_n =\bigoplus_{j>2}(\mathfrak{g}_{e_2-e_j} \oplus
    \mathfrak{g}_{e_2+e_j} ) \oplus \mathfrak{g}_{e_2} \oplus
    \mathfrak{t}_1 \oplus \mathfrak{t}_2.$$ In the latter, it is $$U_n= \bigoplus_{j>2}(\mathfrak{g}_{e_2-e_j} \oplus
    \mathfrak{g}_{e_2+e_j} ) \oplus \mathfrak{t}_1 \oplus
    \mathfrak{t}_2. \phantom{WW}$$

\item[(*)]In this special case,
$
 U_3 = \mathfrak{g}_{e_1-e_2} \oplus
    \mathfrak{g}_{e_1+e_3} \oplus \mathfrak{g}_{e_2}.
$
\end{itemize}

\label{lemma:Un}
\end{lemma}

\begin{proof}
Form a decomposition $\mathfrak{b} = \mathfrak{b}_1 \oplus
\iota(\mathfrak{b}^\downarrow)$ that is compatible with the root
space decomposition.  For $B\in \mathfrak{b},$ write $B=B_1+B_2$
with $B_1 \in  \mathfrak{b}_1$ and $B_2 \in
\iota(\mathfrak{b}^\downarrow)$.  We will write $f_T$ for $f_T^X$.  Note that $B \in
\mathfrak{b}_{f_T}$ if and only if
\begin{equation}\label{eq:ft}
[B,f_{T}]=0.
\end{equation}
To describe $U_n,$ we assume that $B_2$ lies in
$\iota(\mathfrak{b}_{f_{T^\downarrow}}),$ i.e. that
\begin{equation}\label{eq:ft2}
[B_2, \iota(f_{T^\downarrow})]=0.
\end{equation}
  We would like to know what additional conditions on $B$ are necessary to make sure that
   it satisfies (\ref{eq:ft}).  If we write
$$B=\sum_{\alpha \in \Delta^+} c_\alpha E_\alpha + \sum_{i \leq n} c_i T_i,$$ then
(\ref{eq:ft}) imposes linear conditions on the coefficients  in the expansion of
 $B$.  If we choose a representative $\alpha$ or $i$ within each linear condition and
 denote the set of representatives by $P$, then
 $$\mathfrak{b}/\mathfrak{b}_f \simeq \bigoplus_{\alpha \in P} \mathfrak{g}_\alpha
 \oplus \bigoplus_{i \in P} \mathfrak{t}_i.$$
  The natural action of $Q_{f_T}$ has the same determinant on both spaces.
Note that we only need to include representatives for linear
   conditions that do not already arise as conditions for (\ref{eq:ft2}).
    We carry out this plan by describing the set of
    representatives in each of the cases.

Case (C1).  In this case, $f_T= \iota(f_{T^\downarrow}).$ Condition
(\ref{eq:ft}) boils down to
\begin{equation}\label{eq:c1}
[B_1,\iota(f_{T^\downarrow})]=0.
\end{equation}
Write $B_1 = \sum_S c_\alpha E_\alpha +c_1 T_1.$  If we expand the
left hand side of (\ref{eq:c1}) in terms of root space coordinates,
the resulting linear conditions imposed by (\ref{eq:c1}) all take
the form $c_\gamma=0$ for $\gamma$ in some set $\Omega$.  The
quotient $U_n$ then takes the form $\bigoplus_\Omega
\mathfrak{g}_\alpha.$  Deciphering (\ref{eq:c1}) explicitly leads to
the description in the statement of the lemma.

Case (C2).  In this case, $f_T= \iota(f_{T^\downarrow}) +
E_{e_1-e_{\tilde{\phi}(n)}} .$  Equation (\ref{eq:ft}) reduces to
\begin{equation}\label{eq:c2}
[B_1,\iota(f_{T^\downarrow})]+[B_1,E_{e_1-e_{\tilde{\phi}(n)}}]+
[B_2,E_{e_1-e_{\tilde{\phi}(n)}}]=0.
\end{equation}
We can again write $B_1 = \sum_S c_\alpha E_\alpha +c_1 T_1$ and
expand (\ref{eq:c2}) in terms of root space coordinates.   For each linear condition on the coefficients
obtained from  (\ref{eq:c2}), we select as representative the
largest root $\gamma$ such that $c_\gamma$ appears within the linear
equation.  If, however, $c_i$ also appears within a linear
condition,  we select the coefficient $i$ instead.  When we account
for linear conditions that already appear in (\ref{eq:ft2}), we
obtain the description of $U_n$ in the statement of the lemma.

Case (N1).  In this case,  $f_T= \iota(f_{T^\downarrow}) +
E_{e_1+e_{\tilde{k}}} .$  Equation (\ref{eq:ft}) reduces to
\begin{equation}\label{eq:n1}
[B_1,\iota(f_{T^\downarrow})]+[B_1,E_{e_1+e_{\tilde{k}}}]+
  [B_2,E_{e_1+e_{\tilde{k}}}]=0.
\end{equation}
In types $B$ and $D$, the method of case (C2) can be used verbatim,
we only have to account for the different linear conditions imposed
by (\ref{eq:n1}). When $X=C$, we merely have to account for the
different definition of $T^\downarrow$ in this case by letting $B_1=
\sum_S c_\alpha E_\alpha +c_1 T_1 +c_2 T_2$ for the appropriate set
$S$.

Case (N2).  In this case,  $f_T= \iota(f_{T^\downarrow}) +
E_{e_1+e_{2}} .$  Equation (\ref{eq:ft}) reduces to
$[B_1,\iota(f_{T^\downarrow})]+[B_1,E_{e_1+e_{2}}]+
  [B_2,E_{e_1+e_{2}}]=0$
and the method of case (C2) can again be used verbatim to describe
$U_n$.

Case (N3).  In type $B$,  $f_T= \iota(f_{T^\downarrow}) +
E_{e_1-e_{2}} + E_{e_1},$ while in type $D,$ $f_T=
\iota(f_{T^\downarrow}) + E_{e_1-e_{2}} + E_{e_1+e_2}.$ In both
cases, $f_{T^\downarrow} = 0$ and $\mathfrak{b}_{f_{T^\downarrow}} =
\mathfrak{b}^\downarrow$. Hence equation (\ref{eq:ft}) reduces to
$
[B, E_{e_1-e_{2}}] + [B, E_{e_1+e_{2}}] =0$ in type $B$
and
$
 [B, E_{e_1-e_{2}}] + [B, E_{e_1+e_{2}}] =0
$
in type $D$.

Case (*).   In this case,  $f_T= \iota(f_{T^\downarrow}) + E_{e_1}
.$Equation (\ref{eq:ft}) reduces to
$[B_1,\iota(f_{T^\downarrow})]+[B_1,E_{e_1}]+
  [B_2,E_{e_1}]=0$
and the method of case (C2) can again be used verbatim to describe
$U_n$.
\end{proof}

\begin{corollary}
For a standard Young or domino tableau $T,$ $$ \dim U_n = \dim \calV_T - \dim
\calV_{T^\downarrow}.$$ \label{cor:dimension}
\end{corollary}
\vspace{-.3in}
\begin{proof}
We can compute $\dim \calV_T - \dim \calV_{T^\downarrow}$ from the
formula for the dimension of a nilpotent orbit, see \cite{cm}.  Let $[\lambda_1,
\cdots \lambda_p]$ be the dual partition to $shape \, T$. In each of
the cases, $\dim \calV_T - \dim \calV_{T^\downarrow}$ equals
$$
\frac{1}{2}(\dim \calO_{shape \, T} - \dim \calO_{shape \,
T^\downarrow}) = \left\{
\begin{array}{ll}
 \lambda_2+\lambda_3              & \text{ Case (C1)} \\
 \lambda_1                        & \text{ Case (C2) and $X$ = $A$}\\
 \lambda_1-1 + \lambda_3          & \text{ Case (C2) and $X$ = $B$ or $D$}\\
 \lambda_1+1                      & \text{ Case (C2) and $X$ = $C$}\\
 \lambda_1                        & \text{ Case (N1) and $X$ = $B$ or $C$}\\
 \lambda_1-1                      & \text{ Case (N1) and $X$ = $D$}\\
 2\lambda_1-1                     & \text{ Case (N2)} \\
 \lambda_1                        & \text{ Cases (N3) and (*)} \\
 \end{array}
\right.
$$
One can now check these are exactly the dimensions of the
corresponding spaces $U_n$.  We detail the calculation in case (C2) for groups of type $C$.  The
other cases are similar.  Recall from \cite{pietraho:springer} the two types of vertical
dominos that arise within a domino tableau, and denote by $I^-$ and
$I^+$ the set of labels of the dominos of that type that are
contained in the tableau $T$.  From our description of
$U_n$, we obtain:
\begin{align*}
\dim U_n = & \phantom{W} |\left\{j<\phi(n)\, |  \, j \notin N^T \cup (\pi_A (T))^2 \right\}| +
            | \left\{ i \in (\pi_A (T))^2 \, | \, \phi(i)>\phi(n) \right\}| \\
 + & \phantom{W}|\left\{j \, | \, j \neq \phi(i) \text{ for } i \in (\pi_A (T))^2, j
\neq n \right\}|
   +3 \\
 = & \phantom{W}\left\{ j \in I^- \, | \,  j \neq \phi(n), \text{ and if }
  j>\phi(n), \text{ then } j \in Im \, \phi \right\}\\
+ & \phantom{W}\left\{ j
\in I^- \cup N^T \, | \,  j \neq n \right\}
= (|I^-|-1)+(|I^-| + |N^T| -1) + 3 \\
= & \phantom{W}2 |I^-| + |N^T| + 1 = \lambda_1+1,
\end{align*}
as claimed.

\end{proof}
\subsection{The Trace of the Adjoint Action}

Let $\mathfrak{t}_f$ be a maximal torus inside the Lie algebra
$\mathfrak{q}_f$.  Again, we abbreviate $f_T^X$ as $f_T$.  It is easy to check that $\mathfrak{q}_{f_T} \cap \mathfrak{t}$ is a maximal
torus in $\mathfrak{q}_{f_T}$.  The inductive procedure of the
previous section provides a description of the coordinates of
$\mathfrak{t}_f$.  The trace of the adjoint action of
$\mathfrak{t}_f$ on $\mathfrak{q}/\mathfrak{q}_f$ can then be
calculated as a sum of the traces of the actions of the quotient
spaces $U_m$ for $m \leq n$.  In keeping with the inductive philosophy of this
section, we compute this trace on the space $U_n$, separating each
of the inductive cases.

\begin{proposition}
Let $f=f_T$ and write  $a \in \mathfrak{t}$ as $a= \sum_{1
\leq i \leq n} a_i \, \mathfrak{t}_i.$  Then $a$ lies in
$\mathfrak{t}_f$ iff $ \; \sum_{2 \leq i \leq n} a_i \,
\mathfrak{t}_{i}$ lies in the torus
$\iota(\mathfrak{t})_{\iota(f^\downarrow)}$ and additionally
\begin{enumerate}[(i)]
\item $a_1 = a_{\tilde{\phi}(n)}$ in case \textup{(C2)},
\item $a_1 = - a_{\tilde{k}}$ in case \textup{(N1)}, where $\{k,n\}$ is a pair in $ N_1^T$,
\item $a_1 = 0$ in cases \textup{(N2)}, \textup{(N3)}, as well as \textup{(*)}.
\end{enumerate}
\label{proposition:tf}
\end{proposition}

\begin{proof}
This follows immediately from the inductive description of $f_T$.
\end{proof}

\begin{proposition}
Let $p$ be the partition corresponding to the nilpotent orbit passing through $f_T$ and let $\lambda=[\lambda_1,
\lambda_2, \ldots, \lambda_m ]$ be its dual partition.  The trace of the adjoint action of $\mathfrak{t}_f$
  on the quotient $U_n$ is
   listed below:
\begin{itemize}
\item[(C1)]  The trace is $ -\lambda_2 \, a_1 + \sum_{i \in T^2} a_{\tilde{\i}}$ in type $A$ and
$-(\lambda_2+\lambda_3) a_1$ otherwise.

\item[(C2)] The trace is $ - \lambda_1 \, a_1 + \sum_{i \in T^1} a_{\tilde{\i}}$ in type $A$ and
$(- \lambda_1 - c)  a_1$ otherwise, with $c=2$ in type $C$ and $c=-2+\lambda_3$  in types $B$ and $D$ when $N_3^T = \emptyset$.

\item[(N1)]  The trace is $ - (\lambda_1-c)a_1$ in types $B$ and $D$ and zero in type $C$.  The constant $c$ is defined as in \textup{(C2)}.

\item[(N2)] The trace is $0$.

\item[(N3)] The trace is $-\lambda_1 a_1.$

\item[(*)] The trace is $-2a_1 - a_3$.
\end{itemize}
\label{prop:trace}
\end{proposition}

\begin{proof}  We use the description of the quotient $U_n$ in Lemma
  \ref{lemma:Un} together with Proposition \ref{proposition:tf}.  In type $A$,
  determining the trace is simply a matter of reading off the
  coordinates.  We provide the calculations in Case (C1) for the other classical
  types which is only a little more subtle. The other cases are similar.

Case (C1).  By reading off the coordinates, we find that the trace
is
$$ \sum_{(\pi_A(T))^2} (-2a_1 + a_{\tilde{\i}}- a_{\tilde{\phi}(i)})
 +  \sum_{i \in
N^T} (-a_1 - a_{\tilde{\i}}) \; \big(+ (-  a_1 + a_{\tilde{3}})\big) $$
\noindent where the final parenthetical expression appears iff some
sub-tableau $T(m)$ of $T$ lies in case (*).  Applying Proposition \ref{proposition:tf}
reduces the above  to $$-2 \, |(\pi_A (T))^2| \, a_1 - \, |N^T| \, a_1
\; (-1) = - (\lambda_2+\lambda_3) a_1,$$
as claimed.
\end{proof}

For future use, let us define the vector $(c_n, c_{n-1}, \ldots,
c_1)$ by letting $c_i$ equal the number of times the term $a_i$
appears in the expression for the trace of the adjoint action on
$\bigoplus_{i \leq n} U_i$ described by Proposition
\ref{prop:trace}.


\section{Infinitesimal Characters}\label{section:infinitesimal}

Armed with the constructions of the previous section, we are ready
to examine the Graham-Vogan construction.  We restrict our work
to those representations that arise from spherical orbital
varieties, when the corresponding Lagrangian coverings are
quotients of $G$.

Let $\calO$ be a spherical nilpotent orbit of a classical simple Lie group
$G$, fix a Borel subgroup $B$,   and consider the orbital variety
$\calV =\calV_T \subset \calO$ that corresponds to the standard Young or domino tableau $T$
by the parametrization of Theorems \ref{theorem:ovA} and \ref{theorem:ovX}.  Write $Q$ for its stabilizer in
$G$ described by Theorem \ref{theorem:tau}.  We will write $f \in \calV$ for the point $f_T$ in type $A$ as well as the point  $f_T^X$ in the other classical types defined in \S \ref{subsection:basepoint} via the map $\phi$ of Proposition \ref{proposition:basepointA}.
Let $Q_f \subset Q$ be its stabilizer.  Lemma \ref{lemma:dense} implies that the $Q\cdot f$ and $B\cdot f$ are dense in $\calV$.


\subsection{Characters, Weights, and Extensions}\label{subsection:extensions}

A crucial step in the Graham-Vogan construction of the space $V(\calV,\pi)$
relies on the existence of a map $j_\pi$, where we use the notation of \S \ref{section:construction}.  It is by no means clear that such a map exists.  The goal of this section is to describe a condition for its existence which we will use in \S \ref{subsection:VVpi}. Graham and Vogan's
construction examines the character $\alpha$ of $Q_f$
given by the square root of the absolute value of the real
determinant of $Q_f$ acting on the tangent space
$\mathfrak{q}/\mathfrak{q}_f$ of $\calV$ at $f$.   A homomorphism
$j_{\pi}$ exists iff there is a representation $\gamma$ of $Q$ such that $\gamma |_{Q_f} \supset \alpha.$ If $\gamma$ is character, then $j_{\pi}$ will be an isomorphism.

We begin by examining the weight $w_\alpha$ of $\alpha$.  First note that $\alpha$ is a real character.
Recall the vector $(c_n, \ldots , c_1)$ defined at the end of the
previous section.  Splitting the weight of $\alpha$ into
holomophic and anti-holomorphic parts, we obtain:
$$
w_\alpha = \textstyle{(\frac{c_n}{2}, \frac{c_{n-1}}{2}, \ldots,
\frac{c_1}{2}) (\frac{c_n}{2}, \frac{c_{n-1}}{2}, \ldots,
\frac{c_1}{2})}.
$$
We interpret a weight of $Q_f$ as an equivalence class of
weights of $Q$ so the above is just a representative of such an
equivalence class. To answer the existence question for $\gamma$, we
examine its corresponding weights.  Suppose
that $\gamma$ is a real character.  Write the Levi subalgebra
$\mathfrak{l}$ as a sum of reductive parts as $\bigoplus_{i\leq s}
\mathfrak{g}(l_i)$.  A real character $\gamma$ of $L$ takes
the form
\begin{equation}\tag{\dag}
\gamma(A) = \prod_{i\leq s} |\det A_i|^{\alpha_i} = \prod_{i \leq s} (\det
A_i)^{\textstyle{\frac{\alpha_i}{2}}}  \overline{(\det
A_i)}^{\textstyle{\frac{\alpha_i}{2}}}\label{eq:gamma}
\end{equation}
where $A \in L$, $\alpha_i \in \mathbb{R}$ and $A_i$ is the
restriction of $A$ to the $i$th reductive part of $L$.  It has weight
$$
w_\gamma = \textstyle{(\frac{\alpha_n}{2}, \frac{\alpha_{n-1}}{2},
\ldots, \frac{\alpha_1}{2}) (\frac{\alpha_n}{2},
\frac{\alpha_{n-1}}{2}, \ldots, \frac{\alpha_1}{2})}.
$$
We would like to know conditions under which $w_\gamma$ lies in the
same equivalence class of weights of $Q$ as $w_\alpha$.   If we write $w_\gamma = w_\alpha + \epsilon$ for some weight
$\epsilon$, then in the case
of spherical orbits Proposition \ref{proposition:tf} implies that this occurs iff
\begin{itemize}
\item $\epsilon_i+\epsilon_j= 0$ whenever $i= \phi(j)$,

\item $\epsilon_i-\epsilon_j= 0$ whenever $\{i,j\} \in N_1^T$, and

\item $\epsilon_i = 0$ for all $i \notin N^T \cup T^2 \cup \phi(T^2).$
\end{itemize}
Denote the set of weights $w_\gamma$ that satisfy the above
conditions by $HW_r(w_\alpha)$ and write $HW_r^1(w_\alpha)$ for the subset of weights which correspond to a real character $\gamma$ of $Q$.  We would like to understand the relationship between these two sets.  First, let us
define some notation.

For a parabolic subgroup $Q$ of $G$, we  group
the coordinates that correspond to the same reductive part of its
Levi $L$ by setting them off with an additional set of parentheses.
If $$\mathfrak{l} = \bigoplus_{i \leq s} \mathfrak{g}(l_i) \text{\hspace{.2in}
and \hspace{.2in}} \mathfrak{g}(l_j) \cap \mathfrak{t} = \bigoplus_{c_i
  \leq j \leq d_i} \mathbb{C} T_j,$$ then we will write  a weight $a$ as $$a = ((a_n \ a_{n-1} \ \ldots a_{d_1})(a_{c_2} \ \ldots \
a_{d_2}) \ldots (a_{c_k} \ldots a_{d_k}) \ldots (a_{d_l} \ldots
a_1)).$$

\begin{proposition}
A weight $a \in HW_r(w)$ lies in $HW^1_r(w)$ iff all coefficients
corresponding to a given reductive part of the Levi of $Q$  are the
same.  That is, iff  $$a_{c_k} = a_{c_k+1} = \ldots a_{d_k} \text{
for all } 1 \leq k \leq l.$$ \label{proposition:onedim}
\end{proposition}

\vspace{-.3in}
\begin{proof}
Suppose that $a$ satisfies the above hypothesis.  Then a character
of $Q$ with weight $a$ is given by a product of exponents of
absolute values of determinants of the reductive parts of $L$.  The
exponent of the determinant of the part corresponding to $\{c_k,
c_k+1, \ldots, d_k \}$ is given by twice their common value, as per
the description of real characters in (\ref{eq:gamma}).
\end{proof}

Now suppose that $\gamma$ is an arbitrary character of $Q$ that
restricts to the real character $\alpha$ on $Q_f$.  Then $\gamma$
takes the form
$
\gamma = \chi \cdot \gamma'
$
where $\gamma'$ is a real character such that $\gamma' |_{Q_f} =
\alpha$, and $\chi$ is a unitary character such that $\chi |_{Q_f} =
1$.  In particular, this means that $\chi |_{T_f} = 1$. If we write
$A \in T$ as $\sum a_i T_i$, then
$$
\chi(A) = \prod_{i \leq n}
\left(\textstyle{\frac{a_i}{|a_i|}}\right)^{\beta_i} = \prod_{i \leq
n} a_i^{\textstyle{\frac{\beta_i}{2}}}
(\overline{a_i})^{\textstyle{-\frac{\beta_i}{2}}} .
$$
It has weight
$$
w_\chi = \textstyle{(\frac{\beta_n}{2}, \frac{\beta_{n-1}}{2},
\ldots, \frac{\beta_1}{2}) (-\frac{\beta_n}{2},
-\frac{\beta_{n-1}}{2}, \ldots, -\frac{\beta_1}{2})}.
$$
The character $\chi$ restricts to the identity on $Q_f$ iff $w_\chi$
lies in the equivalence class of $0$ of weights of $Q$.  Again by Proposition \ref{proposition:tf}, this occurs
iff
\begin{itemize}
\item $\beta_i+\beta_j= 0$ whenever $i= \phi(j)$,

\item $\beta_i-\beta_j= 0$ whenever $\{i,j\} \in N_1^T$, and

\item $\beta_i = 0$ for all $i \notin N^T \cup T^2 \cup \phi(T^2).$
\end{itemize}
Furthermore, because $\chi$ is a unitary character of $L$,
its entries also  satisfy the conditions of Proposition
\ref{proposition:onedim}.
For a weight $w_\gamma$, write $w_\gamma^h$ and $w_\gamma^a$ for its
holomorphic and anti-holomorphic parts.  Note that  $w_\gamma$ will be in the same equivalence class as $w_\alpha$ iff $w_\gamma^h = (d_n, d_{n-1}, \ldots , d_1)$
satisfies:
\begin{itemize}
\item $d_i+d_j= c_i + c_j$ whenever $i= \phi(j)$,

\item $d_i-d_j= c_i - c_j$ whenever $\{i,j\} \in N_1^T$, and

\item $d_i = c_i$ for all $i \notin N^T \cup T^2 \cup \phi(T^2).$
\end{itemize}
Furthermore, such  $w_\gamma$ will be a weight of a character of $Q$  iff the $d_i$ satisfy the
conditions of Proposition \ref{proposition:onedim}.

\begin{definition} Let $w$ be the weight of a one-dimensional
  representation of $Q_f$ and define $HW(w)$ to be the set of
  weights of representations of $Q$ that restrict to $w$ on the torus
  $\mathfrak{t}_f$.   Furthermore, let $HW^1(w)$  be the
  set of weights in $HW(w)$ that correspond to weights of
  characters of $Q$.
\end{definition}

The arguments of this section reduce the question of extending a character $\alpha$ of $Q_f$
 to a description of the set $HW^1(w_\alpha)$.  Propositions
\ref{prop:trace} and \ref{proposition:onedim}
calculate the weight $w_\alpha$  and a
character $\gamma$ that restricts to $\alpha$ on $Q_f$ exists
whenever $HW^1(w_\alpha)$ is non-empty.


\subsection{The Infinitesimal Characters $IC^1(\calO)$}\label{subsection:characters}

The goal of this section is to describe a set of infinitesimal characters that ought to be attached to an arbitrary nilpotent orbit of a classical simple Lie group.  We follow the work of W.~M.~McGovern \cite{mcgovern}.

A classification of unitary representations of complex reductive Lie
groups can be obtained from a construction that begins with a set of
 special unipotent representations first
suggested by Arthur \cite{barbasch}.  However, only special nilpotent orbits arise as associated varieties of  special unipotent representations.
To remedy this shortfall, \cite{mcgovern} suggests extending this set to a set of $q$-unipotent representations. We first recall McGovern's description of the infinitesimal characters of
$q$-unipotent representations for classical groups.  However,  not
all $q$-unipotent infinitesimal characters obtained by his method
can reasonably correspond to representations attached to nilpotent orbits.  After describing this phenomenon more closely,
 we prune the set of $q$-unipotent infinitesimal characters to a set
 that should be attached to nilpotent orbits.

\subsubsection{Infinitesimal Characters of $q$-unipotent
Representations}

We reproduce the procedure from \cite{mcgovern} for attaching infinitesimal characters to
nilpotent orbits.  Given a nilpotent orbit $\calO,$ we first
describe a way of producing an element $h_{\calO}$ in a Cartan
subalgebra of $\mathfrak{g}.$

\begin{proposition}
For each nilpotent element $f \in \mathfrak{g}$, there is a
homomorphism $\psi: \mathfrak{sl}_2 \longrightarrow \mathfrak{g}$
that maps the element $E_{e_1-e_2}$ onto $f.$  If the nilpotent orbit $\calO$
through $f$ corresponds to the partition $[p_1, \ldots, p_l]$, then
 $ h_{\calO}= \psi (T_1-T_2)$ has eigenvalues:
$$ p_1-1, p_1-3, \ldots, -(p_1 -1), p_2-1, \ldots, p_l-1, \ldots ,
-(p_l-1).$$
We can describe the element $h_{\calO}$ more precisely in
terms of its coordinates.

\begin{enumerate}[(i)]
\item If $\mathfrak{g}$ is of type $A$, then the coordinates of
  $h_{\calO},$ regarded as an element of a Cartan subalgebra of
  $\mathfrak{g}$, are its eigenvalues in non-increasing order.
\item If $\mathfrak{g}$ is of type $B, C$, or $D$, embed it in some
  $\mathfrak{sl}(n)$ via the standard representation.
  \begin{enumerate}[a.]
    \item Suppose the partition of $\calO$ has the numeral \textup{I} or
      none at all.  Also suppose that  $0$ occurs as an eigenvalue of
      the matrix $h_{\calO}$  with multiplicity $k$.  Then the
      coordinates of $h_{\calO}$
      are its positive eigenvalues together with $[k/2]$ zeros,
      arranged in non-increasing order.
    \item If the numeral of $\calO$ is \textup{II}, then the coordinates of
      $h_{\calO}$ are obtained in a similar manner, except that the
      final coordinate is replaced by its negative.
  \end{enumerate}
\end{enumerate}
\end{proposition}

\begin{definition}[\cite{bv}]
An irreducible representation of $G$ is {\it special unipotent} if its
annihilator is of the form  $J_{max}(\lambda_{\calO})$ for $\lambda_{\calO}
=
 \frac{1}{2} h_{\calO}.$
\end{definition}

In each of the classical types except for type $B$, let $n$
be the dimension $d$ of the standard representation of ${}^LG$.  In
type $B$, let $n=d+1.$

\begin{definition}[\cite{mcgovern}]
Let $\mathcal{U}$ be a nilpotent orbit in $\mathfrak{sl}(n)$ and $\lambda_\mathcal{U} = \frac{1}{2}
h_\mathcal{U}$.  Let $\lambda'_\mathcal{U}$ be any $SL(n)$-conjugate
of $\lambda_\mathcal{U}$ lying inside a Cartan subalgebra of ${}^L
\mathfrak{g}.$  When regarded as an infinitesimal character of
$\mathfrak{g},$ $\lambda'_\mathcal{U}$ is called {\it q-unipotent}.
\end{definition}

It remains to attach a nilpotent orbit in $\mathfrak{g}^*$
to each of the $q$-unipotent infinitesimal
characters.  The philosophy of the orbit method dictates that this
is the open orbit $\calO$ contained in the associated variety of
$U(\mathfrak{g})/J_{max}(\lambda'_\mathcal{U}).$  We describe it presently:

\begin{theorem}[\cite{mcgovern}]
Suppose that the orbit $\mathcal{U} \subset \mathfrak{sl}(n)^*$
corresponds  to the partition $p$.  The open orbit $\calO$ in the
associated variety of
$U(\mathfrak{g})/J_{max}(\lambda'_\mathcal{U})$ has partition:
\begin{enumerate}[(i)]

\item $p^t$ in type $A,$

\item $(p^t)_B$ in type $B,$

\item $(l(p^t))_C$ in type $C,$

\item $(p^t)_D$ in type $D,$ except when $p$ is very even, in which case $\calO$ depends
 on the choice of $\lambda_\mathcal{U}$ and can be either $(p^t, I)$ or $(p^t,II).$
\end{enumerate}
\end{theorem}
The maps $p_X$ are the $X$-collapses of the partition $p$ and
$l(p)$ is the partition obtained from $p$ by subtracting $1$ from
its smallest term.  By letting $M(\mathcal{U})=\calO$, we can
define a map
$$M: \text{nilpotent orbits in $\mathfrak{sl}(n)$} \longrightarrow \text{nilpotent orbits in $\mathfrak{g}$}.$$
We will interpret $M$ as a map on
partitions.  We also adopt some notation for
a $q$-unipotent infinitesimal character by associating it with the
partition of the type $A$ orbit that is used to compute it.  For
example, the orbit $\calU_{[4^2,1]}$ lies in the preimage
$M^{-1}(\calO_{[2^4]})$ of the type $C$ orbit with partition $[2^4].$
Then
$\lambda'_\calU = (\textstyle{\frac{3}{2},\frac{3}{2}, \frac{1}{2},\frac{1}{2}})$
which we write as $ \lambda'_\calU =[4^2,1].$  This expression is
unique as long as the classical type of the orbit $\calO$ is specified. In the
case of very even orbit in type $D$, we take this to mean that the
infinitesimal character with all nonnegative terms is attached to the
orbit with numeral I and the infinitesimal character with one
negative term is attached to the orbit with numeral II.

According to our present philosophy, the
$q$-unipotent infinitesimal characters that are attached to the
nilpotent coadjoint orbit $\calO$ for a classical
$\mathfrak{g}$ consist of
$$IC(\calO) = \{ \lambda'_\mathcal{U} \; | \; \mathcal{U} \in M^{-1}(\calO) \}.$$
 We describe this set explicitly for
spherical nilpotent orbits.

\begin{proposition}
Let $\calO =\calO_p$ be a spherical nilpotent orbit in a classical
Lie algebra $\lieg$. In type $A$, $IC(\calO_p)=  \{p^t\}.$  In the other classical types,
the set $IC(\calO_p)$ is as follows:

{\bf Type B}
$$
\begin{array}{|l|l|l|}
\hline
\hspace{.8in} p & \hspace{.8in} IC(\calO_p)\\
\hline
[2^{2k},1^{2n-4k+1}]\; k
     \neq \frac{n}{2} &\{[2n-2k+1,2k], [2n-2k,2k+1]\}\\

      [2^{2k},1^{2n-4k+1}] \text{ $n$ even}  & \{[n+1,n]\}\\

       [3,1^{2n-2}]  & \{[2n-1,1^2],[2n-2,2,1],[2n-2,1^3]\}\\

       [3,2^{2k},1^{2n-4k-2}]\; k \neq \frac{n-1}{2},0  & \{[2(n-k)-1-\epsilon ,2k+1+\epsilon,1] \, | \, \epsilon = 0, \, 1\} \\

       [3,2^{n-1}] \text{ $n$ odd}  & \{[n^2,1]\} \\
       \hline
\end{array}
$$

{\bf Type C}
$$
\begin{array}{|l|l|l|}
\hline
\hspace{.8in} p & \hspace{.8in} IC(\calO_p)\\
\hline
 [1^{2n}]  & \{[2n+1]\}\\

        [2,1^{2n-2}]  & \{[2n,1]\}\\

        [2^k,1^{2n-2k}] \; k
     \neq 1 \text{ or }n  \phantom{WW}& \{[2n-k+1,k],[2n-k+1,k-1,1]\}\phantom{WW}\\

        [2^n]  & \{[n+1,n],[n^2,1],[n+1,n-1,1]\}\\
 \hline
 \end{array}
 $$

{\bf Type D}
$$
\begin{array}{|l|l|l|}
\hline
\hspace{.8in} p & \hspace{.8in} IC(\calO_p)\\
\hline
  [2^{2k},1^{2n-4k}] \; k
     \neq \frac{n}{2}  & \{[2n-2k,2k], [2n-2k-1,2k+1]\}\\

       [2^{n}] \text{ $n$ even}  & \{[n^2]\}\\

      [3,1^{2n-3}]  & \{[2n-2,1^2],[2n-3,2,1],[2n-3,1^3]\}\\

       [3,2^{2k},1^{2n-4k-3}] \; k \neq \frac{n-2}{2} \phantom{W}  & \{[2(n-k-1)-\epsilon ,2k+1+\epsilon ,1] \, | \, \epsilon = 0, 1\}\\

       [3,2^{n-2},1] \text{ $n$ even}  & \{[n,n-1,1]\}.\\
\hline
\end{array}
$$
\end{proposition}
\begin{proof}
The proof is much simpler than the statement.  It consists of
understanding the above map and analyzing all the possibilities.
The details are left to the interested reader.
\end{proof}

Unfortunately, even among this list, there already appear orbits
$\mathcal{U}$ whose associated $q$-unipotent infinitesimal
characters $\lambda'_\mathcal{U}$ cannot be attached to the
nilpotent orbit $\calO \subset \mathfrak{g}^*$ in any reasonable
way. To explain this, write $U$ for the spherical $q$-unipotent
bimodule $U(\lieg)/J_{max}(\lambda'_{\mathcal{U}})$ and define the
 $m_{\lambda'_{\mathcal{U}}}$ to be the multiplicity of the associated variety
$\calV(U)$ in the characteristic cycle $Ch(U)$.
The orbit method dictates that in order for $U$ to correspond to a
cover of a nilpotent coadjoint orbit $\calO$ ,
$m_{\lambda'_{\mathcal{U}}}$ cannot exceed the order of the
fundamental group of $\calO$.  That is, $U$ should not be too large
to meaningfully correspond to $\calO$. It turns out that for certain
$\mathcal{U},$ this unfortunately does occur.  Examples of this
phenomenon arise already among spherical nilpotent orbits.

\begin{example}
Let $\mathcal{U}$ be the nilpotent orbit corresponding to the
partition $[6,3]$ in $\mathfrak{sl}(9).$ Fix the type of the Lie
algebra $\lieg$ to be $C$.  Then $M(\mathcal{U}) = \calO_{[2^3,1^2]}
\subset \mathfrak{sp}(8)^*.$  Furthermore,
$\textstyle{\lambda'_{\calU} =
(\frac{5}{2},\frac{3}{2},1,\frac{1}{2})}.$ However,
$m_{\lambda'_{\mathcal{U}}} = 4$ while $|\pi_1(\calO_{[2^3,1^2]})|
=2$.  According to the above philosophy,
$m_{\lambda'_{\mathcal{U}}}$ should not be the infinitesimal
character  of a unipotent representation attached to
$\calO_{[2^3,1^2]}$.  In
 fact, this is also true for any $\mathcal{U}$ with partition of the form $[2n-k+1, k].$
 \end{example}

There are similar examples in the other classical groups not of type $A$.
Therefore, in order to find the set of the
infinitesimal characters of representations attached to the
nilpotent orbit $\calO$, we have to prune the set $IC(\calO).$

\subsubsection{Pruning of $IC(\calO)$}  We would like
to exclude the infinitesimal characters which arise from those nilpotent orbits for which $m_{\lambda_{\calU}'} > |\pi_1(\calO)|$.   Write a partition $p$ as
$[p_1,p_2, ..., p_l]$ and define
\begin{enumerate}
\item[] $a=$ number of distinct odd $p_i$,

\item[] $b=$ number of distinct even nonzero $p_i$, and

\item[] $c= gcd(p_i)$.
\end{enumerate}

\begin{proposition}[\cite{cm}]
Let $\calO =\calO_p$ be an orbit in a classical simple Lie algebra $\lieg$.
Then the order of the fundamental group $\pi_1{(\calO_p})$ is:

\begin{enumerate}[(i)]

\item $c$ in type $A$,

\item $2^a$  in type $B$ if $p$ is rather odd and $2^{a-1}$ otherwise,

\item $2^b$ in type $C$,

\item $2 \cdot 2^{\max(0,a-1)}$ in type $D$ if $p$ is rather odd and $2^{\max(0,a-1)}$ otherwise.
\end{enumerate}
\end{proposition}

We follow \cite{mcgovern} in determining the multiplicity
$m_{\lambda_{\calU}'}$. The process is a bit complex and requires
notation incompatible with some used here, so rather than referring
the reader to the above, we replicate the relevant parts
here using new notation.

\begin{definition}
Let $M(\calU) = \calO$, and suppose that $\calU$ corresponds to the
partition $p$. In each of the  classical types $X= B,C$, and $D,$ we
define two numbers $\mu$ and $\nu.$
\begin{itemize}
        \item When $X=D,$ let $q=p_{odd}=(q_1^{\lambda_1}, \ldots, q_t^{\lambda_t})$
         and break it up into chunks as follows.  Starting from the left, each chunk takes
         on one of the forms: $(q_i^{\lambda_i}, q_{i+1}^{\lambda_{i+1}})$ with both
          $\lambda_i
$ and $\lambda_{i+1}$ odd; $(q_i^{\lambda_i})$ with $\lambda_i$
even; or $(q_i^{\lambda_i})$ with $\lambda_i$ odd and
$\lambda_{i+1}$ even.  Let $\nu$ be the number of chunks of the
first two types.  The number $\mu$ is defined the same way but with
$q=[ (p_{even})_D]_{odd}.$
        \item When $X=B,$ break up $p_{odd}$ into chunks as in type $D$.  Let $\nu_1$
        be the number of chunks of the first type.  Let $c$ be the leftmost chunk of the
         third type and let $\nu_2$ be the number of chunks of the second type to the right
          of
$c$, plus one.  If no $c$ exists, let $\nu_2=0$.  Finally, let
$\nu=\nu_1+\nu_2.$  The number $\mu$ is defined the same way but
with $([r(p_{even})]_{B})_{odd}$.
        \item When $X=C,$ define $\nu$ in the same way as in type $B.$  To define
        $\mu$, replicate its definition in type $D$ but with the partition
        $[(p_{even})_D]_{odd}$.
\end{itemize}
Finally, in each of the cases let $\nu^*= \max(0,\nu-1)$ and $\mu^*=
\max(0,\mu-1).$
\end{definition}

\begin{definition}
To start, write the coordinates of the infinitesimal character $\lambda'_\calU$
as $ ((\textstyle{\frac{i}{2})^{r_i}}, \ldots , (\frac{1}{2})^{r_1}, 0^{r_0}).$
If $\lambda'_\calU$ contains the coordinate
$\textstyle{-\frac{1}{2}},$ write this as an additional
$\textstyle{\frac{1}{2}}.$  In type $B,$ define the following numbers:
        \begin{itemize}
                \item $\kappa = $ number of even positive $i$ with $r_i$ odd and $r_{i-1}$
                even,
                \item $\kappa_1=$ number of even  positive $i$ with $r_i$ odd, $r_{i-1}$
                even, and
                either $r_{i-2} > r_i$ with $i>2,$ or $r_0 > \textstyle{\frac{1}{2}}r_2,$

                \item $\kappa_2=$ number of even  positive $i$ with $r_i$ odd, $r_{i-1}$
                even positive, and the largest integer $j$ with the following property is
                even: for even $m$, $i\leq m \leq j,$  $r_m$ is odd, while for odd $m$ in
                the same range,
$r_m$ is positive even.
        \end{itemize}
In type $D$, first let $i_0$ be the smallest odd integer $i$ with $r_i$ odd if one
 exists.  Otherwise, let $i_0 = \infty.$ Then define:
        \begin{itemize}
                \item  $\kappa = $ number odd $i$ with $r_i$ odd and either $r_{i-1}$ even
                 or $i=i_0$,
                \item $\kappa_1=$ number of odd $i>i_0$ with $r_i$ odd, $r_{i-1}$ even,
                and either $r_{i-2} > r_i,$
                \item $\kappa_2=$ number of odd $i.i_0$ with $r_i$ odd, $r_{i-1}$ even
                positive, and the largest integer $j$ with the following property is odd:
                for even $m$, $i\leq m \leq j,$  $r_m$ is positive even, while for odd $m$
                in the same range
, $r_m$ is odd.
        \end{itemize}
In type $C,$ the definition is a bit longer.  Define a string of integers $i,\ldots j$ to be {\it relevant} if
 $j>i\geq 0$,
 for $i<m \leq j,$ $ r_m$ is odd,
either $i>0$ and $r_i$ is odd, or $i=0$ and $r_i= \textstyle{\frac{1}{2}}(r_{i+2} -1)$, and
the string is maximal subject to the above. Now let
        \begin{itemize}
                \item $E_S$ be the set of positive even integers $i$ in $S$ such that $r_i>1$ and  $i>2$, or $r_{i-1} \neq 1$,
                \item $F_S $ be the set of odd integers $i$ in $S$ with $r_i >1$, and
                \item $\kappa'_S = \max (|E_S \cup F_S| - (length(S) -2), 0).$
        \end{itemize}
We can now list the relevant strings as $S_1, \ldots , S_r$ in such
a way that the ones with $\kappa'_S =2$ come first, followed by the
ones with $\kappa'_S =1,$ and then the ones with $\kappa'_S =0.$
Enumerate the integers in $\cup_S E_S$ as $i_1, \ldots
 i_s$ in such a way that the ones in  $S_1$ come first, etc.
Now let $\kappa(i_a) = 1$ iff $a \leq \nu^*$ and $0$ otherwise.
Also let $\kappa(j_b) = 1$ iff $b \leq \mu^*$ and $0$ otherwise.
Finally, for each relevant string $S,$ we can define
\begin{itemize}
\item $\kappa_S = \sum_{i_a \in E_S} \kappa(i_a) + \sum_{j_b \in F_S} \kappa(j_b).$
\end{itemize}
\end{definition}

We are now ready to describe the multiplicity
$m_{\lambda'_\calU}.$ Let
 $n_B = 2\kappa - \min(\nu^*,\kappa_1) - \min(\mu^*,\kappa_2),$
 $n_C = \sum_S \max(length(S) -2-\kappa_S,0),$ and
$n_D = 2\kappa - \min(\mu^*, \kappa_1) - \min(\nu^*, \kappa_2) +\kappa_3.$

\begin{proposition}[\cite{mcgovern}]
Consider the type $A$ nilpotent orbit $\calU=\calU_q.$ With notation
as above, $m_{\lambda_{\calU}'}$ equals $1$ in type $A$, $2^{n_B}$ in type $B$, $2^{n_C}$ in type $C$, and $2^{max(n_D-2,0)}$in type $D$.
\end{proposition}

\begin{corollary}
Consider a spherical nilpotent orbit $\calO$ and let $M(\calU_p)
=\calO$.  Then

\begin{enumerate}[(i)]

\item $n_B = 2\kappa$ except when $p=[2n-1,1^2],$ or $[2(n-k)-1,2k+1,1]$,
        in which case it equals $ 2\kappa -1$

\item $n_C = \kappa -1$ when $q$ has the form $[2n-k+1,k]$, and is  $0$ otherwise,

\item $n_D= 2\kappa$.
\end{enumerate}
\end{corollary}

\begin{proof}

In type $B$, both $\mu$ and $\nu$ are less than $2,$ except when
$p=[2n-1,1^2]$, $[2(n-k)-1,2k+1,1]$, $[2n-2,2,1],$ or $[2n-2k-2,2k+2,1]$.
In the case of the former two, $\min(\nu^*, \kappa_1)=1,$ and in the
case of all four, $\min(\mu^*,\kappa_2)=0.$ For spherical orbits of
type $D$, both $\mu$ and $\nu$ are less than $2$.  Finally, in type $C$,
relevant strings of length greater than $2$ arise only when
$p=[2n-k+1,k].$
\end{proof}

We are now ready to state a second approximation to the set of
infinitesimal characters that should appear as infinitesimal
characters of representations attached to spherical nilpotent
coadjoint orbits.  For a given nilpotent orbit
$\calO$ of a given classical type, this is the set of characters of
the form $\textstyle{\lambda'_{\calU}}$ with $M(\calU)=\calO$ that
also satisfy the condition $$\textstyle{m_{\lambda'_\calU} \leq
|\pi_1(\calO)|}.$$ We will denote this set by $IC^1(\calO)$ and
compute it in the next proposition.

\begin{proposition}\label{proposition:IC1}
Let $\calO_p$ be a spherical nilpotent orbit in a classical Lie
algebra $\lieg$ that corresponds to the partition $p$. The set
$IC^1(\calO_p)$ of infinitesimal characters attached to  $\calO_p$
by the above procedure is  $\{p^t\}$ in type $A$ and as in the following tables for the other classical types:

{\bf Type B}
$$
\begin{array}{|l|l|l|}
\hline
\hspace{.8in} p & \hspace{.8in} IC^1(\calO_p)\\
\hline

[2^{2k},1^{2n-4k+1}] & \{[2n-2k,2k+1]\} \\

[3,1^{2n-2}] n \neq 2  & \{[2n-2,2,1],[2n-2,1^3]\}\\

[3,1^{2}]  & \{[2^2,1],[2,1^3],[3,1^2]\}\\

[3,2^{2k},1^{2n-4k-2}]\; k \neq \frac{n-1}{2},0  & \{[2n-2k-2,2k+2,1]\}\\

[3,2^{n-1}]   & \{[n^2,1]\}\\
\hline
\end{array}
$$

{\bf Type C}
$$
\begin{array}{|l|l|l|}
\hline
\hspace{.8in} p & \hspace{.8in} IC^1(\calO_p)\\
\hline

[1^{2n}] & \{[2n+1]\}\\

[2^k,1^{2n-2k}] \; k \neq 2 \phantom{WWWW} & \{[2n-k+1,k-1,1]\}\phantom{WW}\\

[2^2,1^{2n-4}]            & \{[2n-1,1^2],[2n-1,2]\}\\

[2^n] n \neq 2  & \{[n^2,1],[n+1,n-1,1]\};\\

[2^2]           & \{[2^2,1],[3,1^2],[3,2]\};\\
\hline
\end{array}
$$

{\bf Type D}
$$
\begin{array}{|l|l|l|}
\hline
\hspace{.8in} p & \hspace{.8in} IC^1(\calO_p)\\
\hline
[2^{2k},1^{2n-4k}]\; k
     \neq \frac{n}{2} \phantom{WWWw}& \{[2n-2k-1,2k+1]\} \phantom{WW}\\

[2^n]     & \{[n^2]\}\\

[3,1^{2n-3}] &\{[2n-3,2,1],[2n-3,1^3]\}\\

[3,2^{2k},1^{2n-4k-3}]  & \{[2n-2k-3,2k+2,1]\}\\
\hline
\end{array}
$$

\label{proposition:mcgovern}
\end{proposition}


\subsection{Infinitesimal Characters of $V(\calV, \pi)$} \label{subsection:VVpi}

Recall the character $\alpha$, defined as the square root of the
absolute value of the real determinant of the $Q_f$ action on
$\mathfrak{q}/\mathfrak{q}_f$ used to define  $V(\calV, \pi)$.
Suppose that $\alpha$ extends to a character $\gamma$ on $Q$.
According to \S\ref{subsection:extensions}, such an extension exists whenever the set
$HW^1(w_\alpha)$ is not empty.  The first goal of this section is to
decide whether and when this occurs.  This is important as the
construction of $V(\calV, \pi)$ relies on the existence of a bundle
isomorphism $j_{\gamma, \pi}$ defined in \S \ref{section:construction}.  In the setting
of spherical nilpotent orbits,  $j_{\gamma, \pi}$ exists precisely
when there is a character $\gamma$ of the parabolic $Q$ which restricts to $\alpha$ on $Q_f$.

The second goal of the section is to decide how well the
infinitesimal characters of $V(\calV, \pi)$ fit within those that
ought to be attached to the nilpotent orbit $\calO$. Suppose that
the half-density bundle on $G/Q$ is given by the character
$\rho_{G/Q}$, and define $\gamma' = \gamma \otimes \rho^{-1}_{G/Q}.$
The space  $V(\calV, \pi)$   is then a subset of $Ind_Q^G
(\gamma')$. If we write $w_\gamma$ for the character of $\gamma$ and
$\rho$ for the half-sum of the positive roots of $G$, then the
associated infinitesimal character is $\chi_\gamma = w_\gamma +
\rho.$ One expects that $\chi_\gamma$ should be a character attached
to $\calO$ in \S \ref{subsection:characters}, that is, it should lie in the
set $IC^1(\calO)$.
We begin with a short list of examples of what is {\it not} true.

\subsubsection{A Few Examples}

First, we show that it is not always possible to find a character
$\gamma$ of $Q$ that restricts to $\alpha$ on $Q_f$.  This occurs
already in type $A$ for the minimal orbit in rank 5.

\begin{example}
Let $\mathfrak{g}= \mathfrak{gl}_5$ and consider the orbital variety
$\calV_T$ associated to the standard Young tableau
$$
\raisebox{3ex}{$T =$ \; \;}
\begin{tiny}
\begin{tableau}
:.1.4\\
:.2\\
:.3\\
:.5\\
\end{tableau}
\end{tiny}
$$

The basepoint $f= f_T = E_{e_2-e_3}$ constructed in Definition
\ref{definition:basepoint}  has dense $B$-orbit by Lemma
\ref{lemma:dense}. The $\tau$-invariant of $T$ and hence that of
$\calV_T$ can be gleaned from Theorem \ref{theorem:tau} and equals $
\{e_1-e_2, e_3-e_4, e_4-e_5 \} .$  If $Q$ is the parabolic
stabilizing $\calV_T$ and $L$ is its Levi subgroup, the
$\tau$-invariant forces $\mathfrak{l}= \mathfrak{gl}_2 \times
\mathfrak{gl}_3$.  We can now compute the weight of the square root
of the absolute value of the determinant of the $Q_f$ action on
$\mathfrak{q}/\mathfrak{q}_f.$ According to the inductive procedure
of Proposition \ref{prop:trace}, the weight of $\alpha$ is $w_\alpha
= -\mathfrak{t}_1 - \mathfrak{t}_3 + \mathfrak{t}_4 +\mathfrak{t}_5$
which we write as $w_\alpha= ((-1 , 0), (-1 , 1 , 1))$  by
grouping terms that correspond to the same reductive part of the
Levi.  The set of weights that restrict to $w_\alpha$ consists of
the one-parameter family $$HW(w_\alpha) = \{w_\alpha + \epsilon =
((-1 , \epsilon_1) , ( -1 -\epsilon_1 , 1 , 1))\}.$$ The weight
$w_\alpha + \epsilon$ corresponds to a one-dimensional
representation of $Q$ iff $ -1 = \epsilon \text{
and } -1-\epsilon = 1.$  This is not possible, implying that
$HW^1(w_\alpha)= \emptyset.$ \label{example:noextension}
\end{example}

One can reasonably expect that the property  $HW^1(w_\alpha)=
\emptyset$ is preserved by induction on tableau. This too is false.

\begin{example}

Let $\mathfrak{g}=\mathfrak{gl}_6$ and consider
$$
\raisebox{3ex}{ $S =$ \; \;}
\begin{tiny}
\begin{tableau}
:.1.4\\
:.2.6\\
:.3\\
:.5\\
\end{tableau}
\end{tiny}
$$
Then $S$ contains $T$ as a subtableau.  Following the procedure of
the previous example, we find that $ w_\alpha= ((-2) , ( 0 , 0 ) ,
(0 , 1 , 1)) $ which extends to a two-parameter set of weights of
the form
$
HW(w_\alpha) = \{w_\alpha + \epsilon = ((-2- \epsilon_1) ,
(\epsilon_1 , \epsilon_2 ) , ( -\epsilon_2 , 1 , 1))\} .$ The
weight $w_\alpha + \epsilon$ corresponds to a one-dimensional
representation $\gamma$ of $Q$ whenever $\epsilon_1=\epsilon_2 =
-1,$ so that $HW^1(w_\alpha)\neq
\emptyset$.  In fact,
\begin{equation*}
w_\gamma = w_\alpha + (1 , -1 , -1 , 1 , 0 , 0) = (\textstyle{\frac{3}{2} , \frac{1}{2} , \frac{1}{2} ,
-\frac{1}{2} ,  -\frac{1}{2} , -\frac{3}{2}}) \in IC^1(\mathcal{O}_{[2^2,1^2]}).
\end{equation*}
\end{example}

One can also hope that if there does exists a character $\gamma$ that
restricts to $\alpha$, then $\chi_\gamma \in IC^1(\calO).$  Unfortunately,
this also fails.

\begin{example} \label{example:notinIC}

Let $\mathfrak{g}= \mathfrak{gl}_6 $ and consider the tableau
$$
\raisebox{2ex}{T= \;  }
\begin{tiny}
\begin{tableau}
:.1.3\\
:.2.5\\
:.4.6\\
\end{tableau}
\end{tiny}
$$
The Levi of the parabolic stabilizing $\calV_T$ is $ \mathfrak{l} =
\mathfrak{gl}_2 \oplus \mathfrak{gl}_2 \oplus \mathfrak{gl}_2$.
Proposition \ref{prop:trace} implies that $ w_\alpha =
\textstyle{((-\frac{3}{2} , -\frac{3}{2}) , ( \frac{1}{2} ,
-\frac{1}{2}) , (\frac{3}{2} , \frac{3}{2}))}. $ The set of weights
that restrict to $\alpha$ on $\mathfrak{t}_f$ is the three-parameter
family
$$
HW(w_\alpha) =  \left\{w(\epsilon_1, \epsilon_2, \epsilon_3)=
\textstyle{((- \frac{3+\epsilon_1}{2}, -\frac{3 + \epsilon_2}{2}), \
(\frac{1+\epsilon_2}{2} , -\frac{1 + \epsilon_3}{2}) , (\frac{3+
\epsilon_3}{2} , \frac{3+\epsilon_1}{2}))}\right\}.
$$
For $w(\epsilon_1, \epsilon_2, \epsilon_3)$ to lie in
$HW^1(w_\alpha)$, we must have $\epsilon_1=\epsilon_2=\epsilon_3$
and $1+\epsilon_2 = -1-\epsilon_3.$  This forces $\epsilon_1 =
\epsilon_2 = \epsilon_3 = -1$.  Hence
$$
HW^1(w_\alpha) = \{ w(-1, -1, -1) = ((-1 , -1) , (0 , 0) , (1 ,
1))\}
$$
which corresponds to the character of the parabolic $Q$ given by
$$
\gamma \begin{tiny}\left(
\begin{array}{lll}
A_1 & * & * \\
0   & A_2 & * \\
0  & 0 & A_3  \\
\end{array}
\right)\end{tiny} = \left( |A_1|^{-2} |A_3|^{2} \right)^{\frac{1}{2}}
$$
\noindent The infinitesimal character of $Ind_Q^G ( \gamma \otimes
\rho^{-1}_{G/Q})$ is then $
\chi_\gamma
=  \textstyle(\frac{3}{2} , \frac{1}{2} , \frac{1}{2} ,
-\frac{1}{2} , -\frac{1}{2} , -\frac{3}{2}).
$
But $\chi_\gamma$ does not lie in $IC^1(\calO_{[3,3]}) =\{(1 \ 1 \ 0 \ 0  -1
\ -1)\}.$  In fact, $\chi_\gamma \in IC^1(\calO_{[4,2]})!$

\end{example}


\subsubsection{Exhaustion of $IC^1(\calO)$}

We address the question of when it is possible to extend the character $\alpha$ of
$Q_f$ to a character of $Q$, and whether the set of such extensions
for a given orbit provides enough candidates whose associated
infinitesimal characters exhaust $IC^1(\calO)$.   Example
\ref{example:noextension} shows that it is certainly not always
possible find an extension $\gamma$ of $\alpha$ for every orbital
variety $\calV \subset \calO$.  However, there exists at least one
orbital variety within each orbit whose associated $\alpha$  does
admit such an extension.  Furthermore, there exists a sufficient
number of such orbital varieties in $\calO$ to account for all
infinitesimal characters in $IC^1(\calO)$.

\begin{theorem}
Let $\calO$ be a rigid spherical orbit or a model orbit with $n>2$ for a classical simple Lie group of rank $n$.
For every $\chi \in IC^1(\calO)$, there exists an orbital variety
$\calV \subset \calO$ for which $\alpha_\calV$ extends to a character $\gamma$ of Q, and
 $\chi_\gamma = \chi$.
\label{theorem:all}
\end{theorem}

\begin{proof} Consider $\calO$ as above.  It is always possible to
   construct a unique standard tableau $T_\calO$ satisfying the following:
        \begin{enumerate}[(i)]
                \item There exists an integer $k$ such that $\forall i \leq k$,                         $i \in T^1$ and $i \notin T^2$,

                \item $k$ is maximal among all standard tableaux of shape equal to
                the partition corresponding to $\calO$.
        \end{enumerate}
When $\calO_{shape \; T_\calO}$ is a very even orbit in type $D$ with Roman numeral II, define a tableau $T_{II}$
by requiring that $\{n-1,n\} = N_1^{T_{II}}$, $ T^1_{II}$ consist of odd numbers, and $T^2_{II}$ consist of even ones.

The desired orbital variety is $\calV_{T_\calO}$ (or $\calV_{T_{II}}$).
The Levi of the stabilizing subgroup of $\calV_{T_\calO}$ has
exactly two reductive components.  We first examine the case where
the largest element of the partition $p$ corresponding to $\calO$ is
$2$ and the Roman numeral associated to $\calO$, if any, is I.   Let
$[\lambda_1, \lambda_2]$ be the partition dual to $p$.  Then
$
w_\alpha = ((c_1, c_1, \ldots c_1),(c_2, c_2, \ldots c_2))
$
where
$$
(c_1,c_2) = \left\{
\begin{array}{ll}
(-\lambda_1, -\lambda_2) & \text{in type $A$,} \\
(-\lambda_1+2, 0)        & \text{in types $B$ and $D$, and}\\
(-\lambda_1-2, 0)        & \text{in type $C$.}
\end{array}
\right.
$$

The elements of $HW(w_\alpha)$ have the form $w(\epsilon_1,
\epsilon_2, \ldots, \epsilon_s) = ((c_1-\epsilon_1, c_1-\epsilon_2, \ldots c_1 - \epsilon_s), (c_2, \ldots, c_2+\epsilon_2, c_2+\epsilon_1))$ if $\calO$ is rigid, $((\epsilon_1, c_1-\epsilon_2, \ldots c_1 - \epsilon_s), (c_2+\epsilon_s, \ldots, c_2+\epsilon_2, c_2+\epsilon_1))$  if {$p = [2^n]$ and $n$ is odd in type $C$}, and $((c_1-\epsilon_1, c_1-\epsilon_2, \ldots c_1 - \epsilon_s), (c_2+\epsilon_s, \ldots, c_2+\epsilon_2, c_2+\epsilon_1))$ otherwise.

In the first case, $w(\epsilon_1, \epsilon_2, \ldots, \epsilon_s)
\in HW^1(w_\alpha)$ iff $\epsilon_i=0$ for all $i$.  In the third
case,  $w(\epsilon_1, \epsilon_2, \ldots, \epsilon_s) \in
HW^1(w_\alpha)$ iff $\epsilon_i=\epsilon_j$ for all $i$ and $j$.
This produces a one-parameter family of weights that depends on the
common value of the $\epsilon_i = \epsilon$.  In the second case,
$w(\epsilon_1, \epsilon_2, \ldots, \epsilon_s) \in HW^1(w_\alpha)$
iff $\epsilon_i=\epsilon_j$ for $i,j \geq 2$ and $\epsilon_1= c_1 -
\epsilon_2$.  This again yields a one-parameter family of weights
that depends on the common value of the $\epsilon_i = \epsilon$ with
$i \geq 2$. Therefore, an orbital variety always exists for which $\alpha_{\mathcal{V}}$  extends to a character $\gamma$ of $Q$.

It is also easy to check that the $\chi_\gamma$ exhaust $IC^1(\calO)$.   In the first case above, that is, whenever $\calO$ is rigid, $|IC^1(\calO)|=1$ and a comparison
with Proposition \ref{proposition:mcgovern} shows that $\{w(0, 0,
\ldots 0)+\rho\} = IC^1(\calO).$  In the third case above in types $B$
and $D,$ $|IC^1(\calO)|=1$ again and with $\epsilon=0$, $\{w(0, 0,
\ldots 0)+\rho \} =  IC^1(\calO).$  In type $C$ when $n>2$,
$|IC^1(\calO)|=2$. Note that $w(-1, -1, \ldots -1) \neq w(0, 0,
\ldots 0)$ and it is an easy check that $\{w(-1, -1, \ldots
-1)+\rho, w(0, 0, \ldots 0)+\rho\} =  IC^1(\calO).$  The second case
is similar.  Now consider the case when $\calO$ is very even in type
$D$ with numeral II and examine the orbital variety
$\calV_{T_{II}}$.  We find that
$$
w_\alpha = \textstyle{\frac{1}{2}}((-2n+2, -2n+2, -2n+4),(-2n+4,
-2n+6), \ldots, (-4, -2),(-2, 0),(0)).
$$
The elements in $HW(w_\alpha)$ have the form $w(\beta, \epsilon_1,
\ldots, \epsilon_s)= \textstyle{\frac{1}{2}}(-2n+2+\beta,
-2n+2+\beta, -2n+4-\epsilon_1),(-2n+4+\epsilon_1, \ldots
(-2+\epsilon_{s-1}, -\epsilon_s),(\epsilon_s)).$ The set of elements
in $HW^1(w_\alpha)$ is a one-parameter family, consisting of
$w(\epsilon_s)= \textstyle{\frac{1}{2}}((\ldots (-4+\epsilon_s
,-4+\epsilon_s),(-\epsilon_s, -\epsilon_s),(\epsilon_s))$. When $n$
is odd, let $\epsilon_s=-1$ and when $n$ is even, let $\epsilon_s =
3$.  Inductively, it is now easy to show that $$w(\epsilon_s)+\rho =
\textstyle{(\frac{n-1}{2}, \frac{n-1}{2},
  \frac{n-3}{2}, \frac{n-3}{2}, \ldots , 3, 3, 1, -1)}$$
which again accounts for $IC^1(\calO)$.

Now assume that $p$ has largest part $3$ and contains parts of size
$2$.  Then according to Proposition \ref{prop:trace},
$
w_\alpha = ((c_1, c_1, \ldots, c_1-1), (c_2, c_2, \ldots c_2))
$ where $(c_1, c_2) = (-\lambda_1+1, 0)$ in both types $B$ and $D$.  The
 elements in $HW(w_\alpha)$ have the form
$$
w(\epsilon_1, \ldots \epsilon_{s+2}) = ((c_1-\epsilon_1, \ldots
c_1-\epsilon_s, c_1-\epsilon_{s+1}), (c_2+\epsilon_{s+2},
c_2+\epsilon_s, \ldots c_2+\epsilon_1))
$$
According to Proposition \ref{proposition:onedim}, $w(\epsilon_1, \ldots
\epsilon_{s+2}) \in HW^1(w_\alpha) $ iff $\epsilon_i=0$ for all $i
\leq s$ and $s+2$, and $\epsilon_{s+1}=-1$.  Furthermore, $\{ w(0,
0, \ldots, 0, -1, 0) + \rho) \} = IC^1(\calO)$.

Now if $p$ has no parts of size two and $n>2$, then $\calO$ is
neither rigid nor model, but the same result holds.  We find that
$
w_\alpha = ((c_1)(c_2, \ldots c_2))
$
and the  form of the elements in $HW(w_\alpha)$ is
$
w(\epsilon_1,\epsilon_{2}) = ((c_1 - \epsilon_1)(\epsilon_2, 0, 0
\ldots 0))
$
where $c_1 = -\lambda_1$.  Now $w(\epsilon_1,\epsilon_{2})$ lies in
$HW^1(\calO)$ iff $\epsilon_2=0$.  It is an easy check that
$\{w(1,0)+\rho, w(0,0)+\rho\} = IC^1(\calO).$  This finishes the
proof of the theorem.

\end{proof}


\subsubsection{Inclusion in $IC^1(\calO)$}

The phenomenon of Example \ref{example:notinIC} fortunately occurs only among certain
model spherical orbits.  For all other spherical orbits, the
infinitesimal character $\chi_\gamma$, if defined, does indeed lie in $IC^1(\mathcal{O})$.
\begin{theorem}
Let $\calO$ be a rigid, non-model spherical nilpotent orbit for a classical simple Lie group and
consider an orbital variety $\calV \subset \mathcal{O}$ with stabilizer $Q$.  Suppose
that there exists a character $\gamma$ of $Q$ which restricts to the
character $\alpha$ on $Q_f$ defined as the absolute value on the
real determinant of its action on $\mathfrak{q}/\mathfrak{q}_f$.
Then  $\chi_\gamma \in
IC^1(\calO)$. \label{theorem:ic}
\end{theorem}
We begin with an example illustrating our approach.
\begin{example}
Let $\mathfrak{g}= \mathfrak{gl}_7$ and let $\calO_{[4,3]}$ be the
nilpotent orbit corresponding to the partition $[4,3]$.  Consider
the standard Young tableau

$$
\raisebox{3ex}{$T=$ \;  }
\begin{tiny}
\begin{tableau}
:.1.4\\
:.2.5\\
:.3.6\\
:.7\\
\end{tableau}
\end{tiny}
$$

The orbital variety $\calV_T$ has stabilizer $Q$ with Levi $L$ whose
Lie algebra is $\mathfrak{l}=\mathfrak{gl}_4 \oplus
\mathfrak{gl}_3.$ We would like to know that if $\gamma$ is a
character of $Q$ which restricts to $\alpha$ on $Q_f,$ then $w_\gamma
+ \rho$ lies in $IC^1(\calO).$  By Proposition \ref{prop:trace} and
the analysis of Section 5.1, $w_\alpha= \textstyle{(( -\frac{3}{2},
-1, -1, -1),(\frac{3}{2}, \frac{3}{2}, \frac{3}{2}))},$ and
$$HW(w_\alpha) = \{ w(\epsilon_1, \epsilon_2, \epsilon_3) =
\textstyle{((-\frac{3}{2} , -\frac{2+\epsilon_1}{2} ,
  -\frac{2+\epsilon_2}{2} , -\frac{2+\epsilon_3}{2}) , (\frac{3+\epsilon_3}{2} ,
  \frac{3+\epsilon_2}{2} , \frac{3+\epsilon_1}{2}))}.$$  Hence $w(\epsilon_1,
  \epsilon_2, \epsilon_3) \in HW^1(w_\alpha) $ iff $\epsilon_i = 1$
  for all $i$.  Therefore  $w_\gamma =w (1, 1, 1)$ and
\begin{equation*}
w_\gamma+ \rho =  \ \textstyle{(\frac{3}{2} , \frac{1}{2} , -\frac{1}{2} ,
-\frac{3}{2} , 1 , 0 , -1)} \in IC^1(\calO_{[4,3]}),
\end{equation*}
as desired.  Now note that
$ w_\alpha + \rho = w(0,0,0) + \rho = \textstyle{(\frac{3}{2} , 1
, 0 , -1 , \frac{1}{2} , -\frac{1}{2} , -\frac{3}{2})}.$ While
$w_\alpha$ does not correspond to a character of $Q$, $w_\alpha+\rho$ is
nevertheless a permutation of $w_\gamma+ \rho$ and also lies in
$IC^1(\calO_{[4,3]})$.  This observation suggests an approach to our
problem.  We will prove: \label{example:first}
\end{example}

\begin{lemma}\label{lemma:ic}
Suppose that we are in the setting of Theorem \ref{theorem:ic}.
Then there exists a weight $w_\beta$ such that
$w_\beta + \rho$ lies in $IC^1(\calO)$ as well as the Weyl group orbit of $w_\gamma +
  \rho.$
\end{lemma}
\noindent The lemma implies that $w_\gamma+\rho \in IC^1(\calO),$
proving Theorem \ref{theorem:ic}.  Its proof will occupy the remainder of this section.

As in our examples, an element of $HW(w_\alpha)$ can be written as  $w(\epsilon_1, \epsilon_2, \ldots \epsilon_s)= (b_n \ b_{n-1} \
\ldots b_2 \ b_1)$.  Each entry $b_i$ may be a constant or
depend on a single independent variable and at most two entries
can depend on the same independent variable.  We can divide the
entries of $b$ into disjoint maximal strings of entries of the form
$b_l, b_{l-1}, \ldots b_k$ which satisfy:
\begin{itemize}
\item $b_l$ and $b_k$ both depend on the same independent variable, and

\item there is no pair $(l', k')$ such that $b_{l'}$ and $b_{k'}$
  both depend on the same independent variable and
        $l'>l$ and $k'<k.$
\end{itemize}

For such a maximal string, call $I = (k, k+1, \ldots, l)$ a {\it
dependent
  interval} of $b$.  It is an easy consequence of Proposition \ref{proposition:tf}
that if $i$ lies in a dependent interval,  the entry $b_i$ is not
constant.  If $i \in I$ and $b_i$ depends on the variable
$\epsilon_{N_i}$,  we will say that $\epsilon_{N_i}$ {\it
corresponds} to $I$.  Note that each $\epsilon_i$ for $i \leq s$
corresponds to one and only one dependent interval $I$.
If for all $ i \leq \textstyle{\frac{l-k+1}{2}}$ the entries
$b_{l-i}$ and $b_{k-i}$ depend on $\epsilon_{N_i}$ for some
$N_i \leq s$, we will say that $I$ is {\it simple}.

For each simple dependent interval $I=(k \ldots l)$, define a
permutation
$$
\sigma_I = \prod_{i < \textstyle{\frac{l-k}{2}}}(l-i \ \
\textstyle{\frac{k+l}{2}}+i)
$$
 as a product of transpositions;   $\sigma_I$ simply interchanges the first half of the
entries of $I$ with the second half,
preserving the relative order of elements within each.  By
hypothesis, we know that there exists a character $\gamma$ of $Q$
which restricts to $\alpha$ on $Q_f$.  Hence there exists a constant
$c_i$ for each variable $\epsilon_i$ such that  $w(c_1, c_2, \ldots)
\in HW(w_\alpha)$ that equals $w_\gamma$.  If there exists a $c_i
\neq 0$ that corresponds to the dependent interval $I$, we say that
$I$ is {\it non-zero.}

\begin{example} Maintain the setting of Example \ref{example:first}.
There is a unique dependent interval $I=(1, \dots, 6)$ corresponding to the entries
$$
 (b_6 \ b_5 \ b_4 \ b_3 \ b_2 \ b_1) = (-\textstyle{\frac{2 +
 \epsilon_1}{2}}, -\textstyle{\frac{2 +\epsilon_2}{2}} ,
 -\textstyle{\frac{2 + \epsilon_3}{2}}, \textstyle{\frac{3 + \epsilon_3}{2}},
 \textstyle{\frac{ 3+ \epsilon_2}{2}}, \textstyle{\frac{3+
\epsilon_1}{2}} ).
$$
In fact, $I$ is  simple and
$
\sigma_I = (3 \ 6) \ (5 \ 2) \ (4 \ 1).
$
Now note that if we write $w_\alpha+\rho$ as $(c_7 \ c_6 \ldots
c_1)$, then $w_\gamma + \rho = (c_{\sigma(7)} \ c_{\sigma(6)} \ldots
    c_{\sigma(1)}).$  Hence, at least in this case, we have produced a
    method of describing the permutation relating $w_\gamma+\rho$ and
    $w_\alpha+ \rho.$

\end{example}

We are ready to define $w_\beta$ and state a proposition outlining the remainder of the proof of Lemma \ref{lemma:ic}.
\begin{definition}\label{definition:wbeta}
We define $w_\beta$ inductively.  For $\delta= \alpha$ or $\beta,$ let $v_\delta = w_\delta -
\iota(w_\delta^\downarrow).$  Further,
let
$$
v_\beta = v_\alpha + \left\{
\begin{array}{ll}
-(n+2)T_1  &  \text{in case (N1) where } X = C \\
-(n+1)T_1  &  \text{in case (N2) where } X = C \\
\hspace{.3in} T_1       &  \text{in case (N3) where } X = B \text{ or } D, \text{ and}\\
\hspace{.3in}  T_3       &  \text{in case (*) where } X= B.
\end{array}
\right.
$$
\end{definition}

\begin{proposition}\label{proposition:final}
If all non-zero dependent intervals in $HW(w_\alpha)$ are simple, then
    \begin{enumerate}[(i)]
        \item \label{proposition:final:i} $w_\beta + \rho \in IC^1(\calO)$,
        \item \label{proposition:final:ii} $w_\gamma + \rho = \sigma (w_\beta+\rho)$, where $\sigma$ is the
  product of the $\sigma_I$ taken over all non-zero simple dependent
  intervals $I$ and acts by permuting the order of the entries of the weights.

 \hspace{-.4in}Furthermore,

        \item \label{proposition:final:iii} a non-zero non-simple dependent interval cannot exist under the
hypotheses of Theorem \ref{theorem:ic}.
\end{enumerate}
\end{proposition}

We will verify the present proposition using a sequence of four lemmas.  First, note that
Proposition \ref{proposition:tf} implies that $w_\beta \in
HW(w_\alpha)$.
 Because of the hypotheses of Theorem \ref{theorem:ic}, we know that there is a
  weight $w_\gamma  \in HW^1(w_\alpha)$.  If we write a general element of $HW(\alpha)$ as
$w(\epsilon_1, \ldots \epsilon_s)= (b_n, b_{n-1}, \ldots, b_1),$
then there exists constants $c_1, \ldots c_s$ such that $w(c_1,
\ldots c_s) = w_\gamma.$  Because $\calO$ is rigid, there exists at
least one entry $b_p$  which is constant.  Note that it does not
belong to any dependent interval.  We will prove:

\begin{lemma}\label{factA}  If $b_p$ is adjacent to a non-zero non-simple dependent
interval, then there are no constants $c_1, \ldots c_s$ such
that $w(c_1, \ldots c_s) \in HW^1(w_\alpha)$.
\end{lemma}

\begin{lemma}\label{factB} If $I_1=(k_1, \ldots l_1)$ is a non-zero non-simple
dependent interval that is adjacent to a simple dependent interval $I_2=(k_2, \ldots l_2)$,
then there are no constants $c_1, \ldots c_s$ such that $w(c_1,
\ldots c_s) \in HW^1(w_\alpha)$.
\end{lemma}

\begin{lemma}\label{factC}If $b_p$ is adjacent to a non-zero simple dependent
interval $I=(k, \dots l)$, then
                $\sigma_I  \left((w_\beta+\rho)_l \ldots (w_\beta+\rho)_k \right) =  \left((w_\gamma+\rho)_l \ldots (w_\gamma+\rho)_k  \right). $
\end{lemma}

\begin{lemma}\label{factD}If $I_1=(k_1, \ldots l_1)$ is a non-zero simple
  dependent interval that is adjacent to either a simple dependent
  interval or a zero non-simple dependent interval $I_2=(k_2, \ldots l_2)$, then
        $\sigma_I  \left((w_\beta+\rho)_{l_1} \ldots (w_\beta+\rho)_{k_1}  \right) = \left((w_\gamma+\rho)_{l_1} \ldots (w_\gamma+\rho)_{k_1}  \right).$

\end{lemma}

Assuming these, we first prove the proposition.

\V \noindent {\bf Proof of Proposition \ref{proposition:final}(\ref{proposition:final:iii}).} If $w(\epsilon_1, \ldots \epsilon_s)$ contains a non-zero non-simple
dependent interval, it must contain at least one that is adjacent to
either a simple dependent interval or a constant.  Lemmas \ref{factA} and \ref{factB}
then provide a contradiction, proving Proposition \ref{proposition:final}(\ref{proposition:final:iii}).

\hfill $\Box$

\V \noindent {\bf Proof of Proposition \ref{proposition:final}(\ref{proposition:final:ii}).} Proposition \ref{proposition:final}(\ref{proposition:final:iii}) shows that $w(\epsilon_1, \ldots \epsilon_s)$ consists
solely of simple dependent intervals and constants.  For an integer
$i$ that either lies in a zero dependent interval or whose
corresponding entry is a  constant, we know that $(w_\gamma+\rho)_i
= (w_\beta+\rho)_i$.  If, however, $i$ lies in a non-zero dependent
interval, Lemmas \ref{factC} and \ref{factD} show that $(w_\gamma+\rho)_i  =
(w_\beta+\rho)_{\sigma(i)}$, which implies Proposition \ref{proposition:final}(\ref{proposition:final:ii}).

\hfill $\Box$

\V \noindent {\bf Proof of Proposition \ref{proposition:final}(\ref{proposition:final:i}).} We would like to show that
$w_\beta \in IC^1(\calO).$  Let $S= \{IC^1(\calO) -
\iota(IC^1(\calO^\downarrow))\}$ and define $w= w_\beta + \rho -
\iota(w_\beta^\downarrow - \rho^\downarrow)$.  It is easy to verify
the lemma for small $n$.  By induction, it is enough to show that $w
\in S$.  The proof in type $C$ includes all the essential elements of
the general proof, and is particularly easy to state.   We detail
each inductive case.
\begin{itemize}
\item[(C1)]  Proposition \ref{prop:trace} implies that $w =
  (\lambda_1, 0 , \ldots, 0)$. Recall the abbreviated character notation of \S \ref{subsection:characters}.
  Note that $k= \lambda_2$ and that the difference
  $w = [2n-k+1, k-1, 1] - [2n-2-k+1, k-1,1]$ always lies in the one or
  two element set $S$.
\item[(C2)] This time, $w= (\lambda_2, 0, \ldots, 0).$  Again using
  the notation of \S \ref{subsection:characters}, we find that
  $w = [2n-k+1,k-1,1] - [2n-k+1, k-3, 1]$ always lies in $S$.
\item[(N1)] Here, $w = (n-2, 0, \ldots, 0)$.  Using the notation of
  \S \ref{subsection:characters}, $w = [n+1,n-1,1] - [n+1,n-3,1]$, which lies in $S$ by
  Proposition \ref{proposition:IC1}.
\item[(N2)] Here, $w = (n-1, 0, \ldots, 0)$.  Using the notation of
  \S \ref{subsection:characters}, $w = [n^2,1] - [n,n-2,1]$, which lies in $S$ by
  Proposition \ref{proposition:IC1}.
\end{itemize}
This accounts for all the cases that arise in type $C$.  For the other
classical types, the proof requires the same inductive verification
except in one instance.  When the partition
corresponding to $\calO$ has no parts of size $1$, then $w_\beta
\notin IC^1(\calO)$.  This is not a contradition, as $\calO$
is not rigid, but it does complicate the induction step.  If
$\mathcal{W}$ is an orbital variety such that
$\mathcal{W}^\downarrow \subset \calO$, then the associated
$w_\beta$ again lies in $IC^1(\calO)$.  This proves Proposition \ref{proposition:final}(\ref{proposition:final:i}).

\hfill $\Box$

 Finally, we address the lemmas.

\V \noindent {\bf Proof of Lemma \ref{factA}.} Write $I = [k, k+1, \ldots l]$
for the non-zero non-simple dependent interval adjacent to $b_p$,
and further assume that $p=k-1$.  The proof for the other
possibility is symmetric.  We utilize notation suggested by Proposition \ref{proposition:onedim}, separating each interval along the break points of the underlying Levi.
The entries of $I$ must have the form:

$$
b_l, b_{l-1}, \ldots b_{m_1}),(b_{m_1-1}, \ldots, b_{m_2}) \ldots
  (b_{m_q -1} , \ldots b_k.
$$
We examine two possibilities.  Either
$l, k-1 \in T^2 \setminus T^1$ and $k \in T^1 \setminus T^2$, or $\{l,k\} \in N_1^T$.
Consider the first case.  The entries of $I$ must then have the form
$$
a_l - \epsilon, a_{l-1}-\epsilon , \ldots, a_{m_1}-\epsilon),
(a_{m_1-1} + \epsilon, \ldots a_{m_2}-\epsilon), \ldots, (a_{m_q-1}
+ \epsilon, \ldots a_k + \epsilon
$$
for some $\epsilon$ since they must correspond to a weight in
$HW^1(w_\alpha)$.  Because all the entries grouped within
parentheses must equal each other, according to Proposition
\ref{proposition:onedim} this gives us the conditions
$a_k+\epsilon  = a_{k-1}
a_{m_i + 1} + \epsilon   = a_{m_{i+1}} - \epsilon$
for all $i<q,$ which translate to
\begin{equation}\tag{a}
\epsilon = a_{k-1} - a_k
 = \textstyle{\frac{a_{m_{i+1}} - a_{m_i+1}}{2}} \hspace{.2in} \text{ for all $i < q.$ }
\label{equation:epsilon}
\end{equation}

 We would like to show that these conditions are
impossible to satisfy.  Proposition \ref{prop:trace} and Definition
\ref{definition:wbeta} give us a description of each of the $a_i$.  We restrict the
proof to type $C$, which contains all the elements of the general
proof.

Let $[\lambda_1(i), \lambda_2(i)]$ be the partition dual to $shape
\, T(i)$. Proposition \ref{prop:trace} implies that
$a_{k-1}  = -\lambda_1(k-1)+2,$
$a_k   = -\lambda_2(k-1) +2,$
$a_{m_2}  = - \lambda_1(m_2)+2,$ and
$a_{m_1+1} = -\lambda_2(m_1+1)$.
Equations (\ref{equation:epsilon}) translate to
\begin{equation}\tag{b}
\epsilon  = -\lambda_1(k-1) + \lambda_2(k-1) \\
 = \textstyle{\frac{-\lambda_1(m_2)+2 + \lambda_2(m_1+1)}{2}}.
\label{equation:epsilon2}
\end{equation}
However, $\lambda_1(k-1) - \lambda_2(k-1) = \lambda_1(l) -
\lambda_2(l)$ because $I$ is a dependent interval.  Furthermore, the
form of the entries in $I$ implies that $\lambda_2(l) >
\lambda_2(m_1+1)$ and $\lambda_1(l)< \lambda_1(m_2)$.  But this
implies that it is impossible to satisfy (\ref{equation:epsilon2}) and we cannot find constants $c_i$ so
that $w(c_1, \ldots c_s) \in HW^1(w_\alpha)$.  The only difference
in proof for the other classical types are the precise values for
the $a_i$.

Now suppose we are in the second case and that $\{k,l\} \in N_1^T$.  The entries corresponding to
the interval $I$ must have the form
$$
a_l + \beta, a_{l-1}-\epsilon , \ldots, a_{m_1}-\epsilon),
(a_{m_1-1} + \epsilon, \ldots a_{m_2}-\epsilon), \ldots, (a_{m_q-1}
+ \epsilon, \ldots a_{k+1} + \epsilon, a_k + \beta.
$$
Because $I$ is non-simple, this means that the interval $\{k+1,
\ldots , l-1\}$ cannot be simple either.  This time, we need to
solve the equations
\begin{equation*}
\epsilon  = a_k + \beta - a_{k+1}
          = a_{l-1} - a_l - \beta
          = \textstyle{\frac{a_{m_{i+1}} - a_{m_i+1}}{2}}
\end{equation*}
First, we find that $\beta =
\textstyle{\frac{(a_{l-1}-a_l)+(a_{k+1}-a_k)}{2}}.$  This means that
we still need to solve
\begin{equation}\tag{c}
\epsilon  = \textstyle{\frac{(a_{l-1}-a_l)-(a_{k+1}-a_k)}{2}}
          = \textstyle{\frac{a_{m_{i+1}} - a_{m_i+1}}{2}}
\label{equation:epsilon3}
\end{equation}
By an analysis similar to the above divided into each classical
type, (\ref{equation:epsilon3}) again cannot be satisfied and Lemma \ref{factA}
holds.

\hfill $\Box$

\V \noindent {\bf Proof of Lemma \ref{factB}.} If $I_2$ is a {\it zero}
interval, then the proof is identical to the proof of Lemma \ref{factA}, as the
only property we needed was the expression for the term $a_{k_1-1}$,
which is the same in the zero case.  Now assume that $I_1$ is to the
left of $I_2$ in the coordinate expression for $w_\gamma$ of this
section; the other possibility has a symmetric proof.  There are
again two cases in the proof.  First assume that $\{k,n\} \notin
N_1^T$.  The two intervals must then have the form
$$
a_{l_1} - \epsilon, a_{l_1-1}-\epsilon , \ldots, a_{m_1}-\epsilon),
(a_{m_1-1} + \epsilon, \ldots a_{m_2}-\epsilon), \ldots, (a_{m_q-1}
+ \epsilon, \ldots a_{k_2} + \epsilon
$$
$$
\text{and \hspace{.5in}}
a_{l_2} -\mu, a_{l_2-1} -\mu, \ldots , a_{m'}-\mu ),(a_{m'-1}+\mu,
\ldots , a_{k_2}+\mu
$$
with the additional restriction that $a_{l_2} - \mu =
a_{k_1}+\epsilon$. Write $\rho$ in coordinates as $(\rho_n ,
\rho_{n-1}, \ldots, \rho_1)$. The proof of Lemmas \ref{factC} and \ref{factD} imply
that either $\mu = 0$, or  $\mu = a_{l_2}-a_{m'} + \rho_{l_2} -
\rho_{m'}$. The first possibility was considered above.  As for the
second, following the outline of the proof of Lemma \ref{factA}, we would like
to solve the equations
\begin{equation}\tag{d}
\epsilon  = a_{l_2} - \mu - a_k
 = \textstyle{\frac{a_{m_{i+1}} - a_{m_i+1}}{2}}  \hspace{.2in} \text{ for all $i<q$.}
\label{equation:epsilon4}
\end{equation}

  In each of the classical types, Proposition
\ref{prop:trace} gives us values for the $a_i$, and we can similarly
give an explicit description of $\rho$.  In a manner similar to the
proof of Lemma \ref{factA}, we can now show that a solution to
(\ref{equation:epsilon4}) does not exists.  A similar analysis works
for the case when $\{k,l\} \in N_1^T$ and Lemma \ref{factB} holds.

\hfill $\Box$

\V \noindent {\bf Proof of Lemma \ref{factC}.}
 Assume that $b_p = b_{k-1}$ as the proof for the
other possibility in symmetric. The entries of $I$ must have the form
$$
b_l, b_{l-1}, \ldots b_m),(b_{m-1}, \ldots b_k.
$$
As in the proof of Lemma \ref{factA}, there are two possibilities.    Either
$l, k-1 \in T^2 \setminus T^1$ and $k \in T^1 \setminus T^2$, or
$\{l,k\} \in N_1^T$.  We examine the first case.  The
second in analogous. Write $\rho$ in coordinates as $(\rho_n ,
\ldots, \rho_1)$. The entries of $w_\gamma$ have the form
$$
a_l -\epsilon , a_{l-1} -\epsilon, \ldots , a_m-\epsilon
),(a_{m-1}+\epsilon, \ldots , a_k+\epsilon
$$
where entries grouped by parentheses must equal since $w_\gamma \in
HW^1(w_\alpha)$. This condition further forces $a_{k-1} = a_k +
\epsilon$, or in other words,
\begin{equation}\tag{e}
\epsilon = a_k - a_{k-1} \label{equation:ep1}
\end{equation}
After examining the definition of the permutation $\sigma_I$, we
need to verify that $
a_{l+i} - \epsilon + \rho_{l+1} = a_{m+i} + \rho_{m+i}$
holds for all $i < (l-k)/2$, which will imply Lemma \ref{factC}.  We first consider
type $A$. First of all, $\rho_{l+i}= n+1 -2(l+i)$, hence we
would like to know whether the equality
$
a_{l_i}-\epsilon + n+1 -2(l+i) = a_{m+i} + n+1 -2(m+i)
$
holds.  Proposition \ref{prop:trace} implies that $a_{l+i}=a_l$ and
$a_{m+i}=a_k$ for all of the above $i$ and the above equation becomes
$a_{l}-a_{k}+{k-l+1} = \epsilon$.  This is possible iff this
equation is compatible with (\ref{equation:ep1}).  To verify this, we
note that repeated application of Proposition \ref{prop:trace}
implies $a_{k-1} = -\lambda_1(k) + (\lambda_1 - \lambda_1(k)) =
\lambda_1 -2\lambda_1(k)$ which also equals $a_l+(k-l+1)$.  This
implies that $a_l -a_k +(k-l+1) = a_{k-1}-a_k$, and thus
Lemma \ref{factC} holds in type $A$.  The proof for groups of other types is
analogous, only complicated by the appearance of horizontal dominos.
However, dominos falling in cases (N2) or (N3) do not affect the
dependent intervals because of Proposition \ref{proposition:tf}.  Case (N1) is
dealt with precisely as in the proof of Lemma \ref{factA}.

\hfill $\Box$

\V \noindent {\bf Proof of Lemma \ref{factD}.} If $I_2$ is a {\it zero}
dependent interval, then the proof is identical to the proof of Lemma \ref{factC}.  We would like to show that in fact, if $I_1$ is a non-zero
simple dependent interval, then $I_2$ must be a zero dependent
interval.  We can assume that $I_1$ is to the left of $I_2$ in the
coordinate notation we have grown accustomed to.  As in Lemma \ref{factC}, the
interval $I_1$ has the form
$$
a_{l_1} -\epsilon , a_{l_1-1} -\epsilon, \ldots , a_m-\epsilon
),(a_{m-1}+\epsilon, \ldots , a_{k_1}+\epsilon
$$
while the interval $I_2$ has the form
$$
a_{l_2} -\mu, a_{l_2-1} -\mu, \ldots , a_{m'}-\mu ),(a_{m'-1}+\mu,
\ldots , a_{k_2}+\mu
$$
with the additional constraint that $l_2 -1 = k_2.$  We would like
to show that $\mu = 0$.  Because $w_\gamma \in HW(w_\alpha)$, we
know that $a_{l_2}-\mu = a_{k_l}+\epsilon.$  But our proof of Lemma \ref{factC}
  implies that in fact,  $a_{l_2} = a_{k_l}+\epsilon,$ forcing $\mu$ to
    be zero, implying Lemma \ref{factD}.

\hfill $\Box$

\end{document}